\documentclass[12pt]{amsart}
\usepackage{amssymb}

\usepackage{tikz}
\tikzset{
  symbol/.style={
    draw=none,
    every to/.append style={
      edge node={node [sloped, allow upside down, auto=false]{$#1$}}}
  }
}

\usepackage[pdfencoding=auto,psdextra]{hyperref}
\usepackage{textgreek}
\usepackage[hmargin=1in,vmargin=1.25in]{geometry}
\definecolor{amethyst}{rgb}{0.6, 0.4, 0.8}
\definecolor{kellygreen}{rgb}{0, 0.4, 0}
\definecolor{americanrose}{rgb}{1.0, 0.01, 0.24}
\definecolor{teal}{rgb}{0, 0.5, 0.5}
\definecolor{lime}{rgb}{0.75, 1, 0}
\definecolor{darklime}{rgb}{0.25, .33, 0}

    \title{A raising operator formula for Macdonald polynomials}

   \author{J. Blasiak}
   \author{M. Haiman}
   \author{J. Morse}
   \author{A. Pun}
   \author{G. H. Seelinger}

    \address[Blasiak]{
    Dept.\ of Mathematics\\
    Drexel University\\
    Philadelphia, PA}
    \email{jblasiak@gmail.com}

   \address[Haiman]{Dept.\ of Mathematics\\
            University of California\\
            Berkeley, CA}
   \email{mhaiman@math.berkeley.edu}

   \address[Morse]{
   Dept.\ of Mathematics\\
            University of Virginia\\
            Charlottesville, VA}
   \email{morsej@virginia.edu}

   \address[Pun] {Dept.\ of Mathematics\\
           Baruch College (CUNY)\\
            New York, NY}
   \email{anna.pun@baruch.cuny.edu}

   \address[Seelinger]{Dept.\ of Mathematics\\
   University of Michigan\\
   Ann Arbor, MI}
   \email{ghseeli@umich.edu}

   \thanks{Authors were supported by NSF Grants DMS-1855784 (J.~B.), DMS-1855804 (J.~M., A.~P. and G.~S.), and DMS-2303175 (G.~S.), and Simons Foundation-821999 (J.~M.).}

   \date{\today}

\newtheorem{thm}{Theorem}[subsection]
\newtheorem{lemma}[thm]{Lemma}
\newtheorem{prop}[thm]{Proposition}
\newtheorem{cor}[thm]{Corollary}
\newtheorem{conj}[thm]{Conjecture}

\theoremstyle{definition}
\newtheorem{defn}[thm]{Definition}

\theoremstyle{remark}
\newtheorem{example}[thm]{Example}

\newtheorem{remark}[thm]{Remark}

\newcommand{\NN}{{\mathbb N}}
\newcommand{\QQ}{{\mathbb Q}}

\newcommand{\ZZ}{{\mathbb Z}}

\newcommand{\kk}{{\mathbf k}}

\newcommand{\bb}{{\mathbf b}}

\newcommand{\Hbold}{{\mathbf H}}
\newcommand{\Kbold}{{\mathbf K}}

\newcommand{\sigmabold}{{\boldsymbol \sigma }}
\newcommand{\nubold}{{\boldsymbol \nu }}

\newcommand{\Ecal}{{\mathcal E}}

\newcommand{\Gcal}{{\mathcal G}}

\newcommand{\ctild}{\tilde{c}}
\newcommand{\Htild}{\tilde{H}}
\newcommand{\Ktild}{\tilde{K}}

\DeclareMathOperator{\inv}{inv}
\DeclareMathOperator{\pol}{pol}

\DeclareMathOperator{\GL}{GL}

\DeclareMathOperator{\SSYT}{SSYT}

\newcommand{\defeq}{\overset{\text{{\em def}}}{=}}

\newcommand{\south}{{\mathrm{south}}}

\newcommand{\leg}{{\mathrm{leg}}}
\newcommand{\arm}{{\mathrm{arm}}}

\newcommand{\HS}{\mathbf{H}}

\newcommand{\xx}{{\mathbf x}}

\newcommand{\zz}{{\mathbf z}}
\newcommand{\dgrm}[1]{{#1}^{*}}
\newcommand{\rdless}{\prec}
\newcommand{\rdlesseq}{\preceq}
\newcommand{\spa}{\hspace{.3mm}}

\newcommand{\str}{{\mathsf{s}}}
\newcommand{\bx}[2]{{\boldsymbol {#1}[#2]}}
\renewcommand{\SS}{\ensuremath{\mathcal{S}}}

\newcommand{\crc}[1]{-#1}

\newlength{\mycellsize}
\newcommand\mytbl[1]{
\vcenter{
\let\\=\cr
\baselineskip=-16000pt \lineskiplimit=16000pt \lineskip=0pt
\halign{&\mytblcell{##}\cr#1\crcr}}}

\newcommand{\mytblcell}[1]{{%
\def \arg{#1}\def \void{}%
\ifx \void \arg
\vbox to \mycellsize{\vfil \hrule width \mycellsize height 0pt}%
\else \unitlength=\mycellsize
\begin{picture}(1,1)
\put(0,0){\makebox(1,1){$#1\vphantom{\crc{#1}}$}}
\put(0,0){\line(1,0){1}}
\put(0,1){\line(1,0){1}}
\put(0,0){\line(0,1){1}}
\put(1,0){\line(0,1){1}}
\end{picture}%
\fi}}

\newlength{\cellsize}
\newcommand\mytableau[1]{
\vcenter{
\let\\=\cr
\baselineskip=-16000pt \lineskiplimit=16000pt \lineskip=0pt
\halign{&\mytableaucell{##}\cr#1\crcr}}}

\newcommand{\mytableaucell}[1]{{%
\def \arg{#1}\def \void{}%
\ifx \void \arg
\vbox to \cellsize{\vfil \hrule width \cellsize height 0pt}%
\else \unitlength=\cellsize
\begin{picture}(1,1)
\put(0,0){\makebox(1,1){$#1\vphantom{\crc{#1}}$}}
\put(0,0){\line(1,0){1}}
\put(0,1){\line(1,0){1}}
\put(0,0){\line(0,1){1}}
\put(1,0){\line(0,1){1}}
\end{picture}%
\fi}}

\DeclareMathSymbol{\shortminus}{\mathbin}{AMSa}{"39}

\newcounter{rownum}
\newcommand{\drawDg}[1]{
      \setcounter{rownum}{0}
      \def\b{#1};
      \foreach \c in \b {
        \draw[thick] (\therownum,0) grid (\therownum+1, \c);
        \addtocounter{rownum}{1};
      }
    }
\newcounter{boxnum}
\newcommand{\drawskewdg}[2]{
  \def\ptn{#1}
  \def\offset{#2}
    \setcounter{rownum}{1}
    \foreach \xstart\xend in \ptn {
      \draw[thick,fill=white] (\xstart+\offset,\therownum-1+\offset)
      grid (\xend+\offset,\therownum+\offset) rectangle (\xstart+\offset,\therownum-1+\offset);
      \addtocounter{rownum}{1}
    }
}

\begin{document}


\begin{abstract}
We give an explicit raising operator formula for the modified
Macdonald polynomials $\Htild_{\mu }(X;q,t)$, which follows from our
recent formula for $\nabla$ on an LLT polynomial and the
Haglund-Haiman-Loehr formula expressing modified Macdonald polynomials
as sums of LLT polynomials.  Our method just
as easily yields a formula for a family of symmetric functions $\Htild^{1,n}(X;q,t)$
that we call $1,n$-Macdonald polynomials, which reduce to a scalar multiple of $\Htild_{\mu}(X;q,t)$ when $n=1$.
 We conjecture that the coefficients of $1,n$-Macdonald polynomials in terms of Schur functions
belong to $\NN[q,t]$, generalizing Macdonald positivity.
\end{abstract}

\maketitle

\section{Introduction} \setcounter{subsection}{1} \label{s intro}
Macdonald polynomials are a family of two-parameter symmetric
functions which have played an important role in developments in
algebra, combinatorics, and geometry since their introduction in the
1980's.  They provide a common generalization of several previously
known families of one-parameter symmetric functions.  Prominent among
these are Hall-Littlewood polynomials, which are classically defined
in terms of raising operators, as in Macdonald's exposition
\cite{Macdonald95}.  Considering that raising operator formulas are
the starting point for the development of Hall-Littlewood polynomials,
the lack of an analogous formula for Macdonald polynomials has been a
hole in the heart of the theory.

Here, we establish a raising operator formula for the modified
Macdonald polynomials $\Htild_{\mu }(X;q,t)$ from which the formula
for Hall-Littlewood polynomials is easily recovered at $q=0$.  Using
notation defined in \S \ref{ss formula for partitions}, our formula
reads
\begin{align} \label{et:raising op main preview}
\omega \Htild_\mu (X;q,t) = \frac{\prod_{\alpha_{ij} \in R_\mu
\setminus \widehat{R}_\mu } \big(1- q^{\arm(\bx{\mu}{i})+1}\spa
t^{-\leg(\bx{\mu}{i})} \mathbf{R}_{ij} \big) \prod_{\alpha_{ij} \in
\widehat{R}_{\mu}} \big(1-q\spa t\, \mathbf{R}_{ij} \big)}
{\prod_{\alpha_{ij} \in R_+} \big(1-q \, \mathbf{R}_{ij}\big)
\prod_{\alpha_{ij} \in R_\mu} \big(1-t \, \mathbf{R}_{ij}\big)} \cdot
s_{1^l}\,.
\end{align}
The proof begins with the Haglund-Haiman-Loehr formula \cite{HaHaLo05}
for $\Htild_{\mu }(X;q,t)$ as a weighted sum of LLT polynomials.  We
then apply the operator $\nabla$, which has $\Htild_{\mu }(X;q,t)$ as
an eigenfunction, and use the formula for $\nabla$ on an LLT
polynomial established in our recent work~\cite{BHMPS-llt}.

  Replacing $s_{1^l}$ in \eqref{et:raising op main preview} with $s_{n^l}$ gives rise to a family of symmetric functions $\Htild^{1,n}_\mu(X;q,t)$, which we call $1,n$-Macdonald polynomials.
We conjecture that the coefficients of these polynomials in terms of Schur functions belong to $\NN[q,t]$, generalizing  Macdonald positivity.
As we will see later,
this can be formulated for all $n$ simultaneously as the statement that a natural
infinite series $\HS_\mu$ of
$\GL_l$ characters has coefficients in  $\NN[q,t]$.

We also
derive a new raising operator formula for the integral form Macdonald
polynomials $J_\mu(X;q,t)$.

Other  raising operator formulas for
Macdonald polynomials have previously appeared in the literature.
Lassalle-Schlosser~\cite{LassalleSchlosser06} inverted the Pieri
formula for Macdonald polynomials $Q_\mu(X;q,t)$ (which differ from
$J_\mu(X;q,t)$ by a scalar factor) to obtain a formula for
$Q_\mu(X;q,t)$ that can be interpreted as a raising operator formula.
Shiraishi~\cite{Shiraishi05} conjectured a similar raising operator
formula for $Q_\mu(X;q,t)$, later proven by Noumi and
Shiraishi in their work \cite{NoumShir12} on the bispectral problem of
the Macdonald-Ruijnesaars $q$-difference operators.  However, these formulas are quite different and more intricate than ours.

\section{Background}
\label{s:background}

\subsection{Partitions and symmetric functions}
\label{ss:background}

The (French style) diagram of a partition $\mu = (\mu_1 \ge \cdots \ge
\mu_k > 0)$ is the set
$\{(i,j) \in \ZZ_{+}^2  :  1\leq j\leq
k,\;1\leq i \leq \mu _{j}\}$.
We identify $(i,j)$ with the lattice
square or \emph{box} whose northeast corner has coordinates $(x,y) =
(i,j)$, and refer to this box as being in {\em column} $i$ and {\em
row} $j$.
We set $|\mu| = \mu _{1}+\cdots +\mu _{k}$ and let $\ell(\mu )=k$ be the
number of non-zero parts.
We write $\mu^*$ for the transpose of a partition $\mu$.
The {\em arm} and {\em leg} of a box $b\in \mu$ are the number of
boxes in $\mu$ strictly east of $b$ and strictly north of $b$,
respectively.

Let $\Lambda = \Lambda(X)$ be the algebra of symmetric functions in
infinitely many variables $X = x_{1},x_{2},\ldots$, with coefficients
in the field $\kk = \QQ (q,t)$.  We follow Macdonald's notation
\cite{Macdonald95} for the graded bases of $\Lambda $, and the
automorphism $\omega \colon \Lambda \rightarrow \Lambda $
given on Schur functions by
$\omega s_{\lambda } = s_{\lambda^*}$.  We also work with series and
symmetric functions in finitely many variables $\zz = z_1,\dots,z_l$.
If $f(X)\in \Lambda$ is a formal symmetric function, then $f(\zz )$ or
$f(z_1, \ldots, z_l)$ denotes its specialization with $X = z_1,
\ldots, z_l,0,0,\ldots$.

Given a symmetric function $f\in \Lambda$ and any expression $A$
involving indeterminates, the plethystic evaluation $f[A]$ is defined
by writing $f$ as a polynomial in the power-sums $p_{k}$ and
evaluating with $p_{k}\mapsto p_{k}[A]$, where $p_{k}[A]$ is the
result of substituting $a^{k}$ for every indeterminate $a$ occurring
in $A$.  The variables $q, t$ from our ground field $\kk $ count as
indeterminates.

By convention, the name of an alphabet $X = x_{1},x_{2},\ldots$ stands
for $x_{1}+x_{2}+\cdots $ inside a plethystic evaluation.  Then $f[X]
= f[x_{1}+x_{2}+\cdots ] = f(x_{1},x_{2},\ldots) = f(X)$.  For
example, the evaluation $f[X/(1-t^{-1})]$ is the image of $f(X)$ under
the $\kk $-algebra automorphism of $\Lambda $ that sends $p_{k}$ to
$p_{k}/(1-t^{-k})$.

The {\em modified Macdonald polynomials} $\Htild_\mu = \Htild _{\mu
}(X;q,t)$ of \cite{GarsHaim96} are defined in terms of the Macdonald
polynomials $Q_{\mu }(X;q,t)$ \cite[VI (4.12)]{Macdonald95} or their
integral forms $J_{\mu }(X;q,t)$ \cite[VI (8.3)]{Macdonald95} by
\begin{equation}\label{e:H-tilde}
\Htild _{\mu }(X;q,t) = t^{\mathsf{n}(\mu )}
J_{\mu}[\frac{X}{1-t^{-1}};q,t^{-1}] = t^{\mathsf{n}(\mu)} \bigl( \prod_{b \in
\mu}(1-q^{\arm(b)+1}t^{-\leg(b)}) \bigr)
Q_{\mu}[\frac{X}{1-t^{-1}};q,t^{-1}],
\end{equation}
where $\mathsf{n}(\mu ) = \sum _{i} (i-1)\mu_{i}$.  The $\tilde
H_\mu(X;q,t)$ also have a direct combinatorial
description~\cite{HaHaLo05}, which we will recall in Theorem~\ref{t
HHL formula}.

When $q=0$, the modified Macdonald polynomials reduce to the
\emph{modified Hall-Littlewood polynomials}
\begin{equation}
 \Htild_\mu(X;0,t) = t^{\mathsf{n}(\mu )} Q_{\mu}[\frac{X}{1-t^{-1}};t^{-1}],
\end{equation}
where the Hall-Littlewood polynomials $Q_{\mu}(X;t)$ are as defined in
\cite[III (2.11)]{Macdonald95}.  At $t=1$ and $t = \infty $, the
$\Htild_\mu(X;0,t)$ specialize to the complete homogeneous symmetric
functions $\Htild_\mu(X;0,1) = h_\mu(X)$ and Schur functions
$\Htild_\mu(X;0,\infty)= s_\mu(X)$.

\subsection{Weyl symmetrization and related operators}
\label{ss:Weyl-etc}

The \emph{Weyl symmetrization operator} $\sigmabold$ for $\GL
_{l}$ is defined by
\begin{equation}\label{e:Weyl-symmetrizer}
\sigmabold (f(z_1,\dots,z_l))
= \sum _{w\in \SS_{l}} w \bigg(
\frac{f(z_1,\dots,z_l)}{\prod _{i < j} (1 - z_j/z_i)}
\bigg)
= \sum _{w\in \SS_{l}} w \bigg(
\frac{f(z_1,\dots,z_l)}{\prod _{\alpha_{ij} \in R_{+}} (1 - z_j/z_i)}
\bigg)
\,,
\end{equation}
where $f \in \kk[z_1^{\pm1},\dots, z_l^{\pm1}]$ is a Laurent
polynomial, $\SS_l$ acts by permuting the variables $z_1, \dots, z_l$,
and $R_+= R_+(\GL_l) = \{\alpha _{ij}  :  1\leq i < j\leq l \}$
denotes the set of positive roots for $\GL _{l}$, with $\alpha_{ij} =
\epsilon_i -\epsilon_j \in \ZZ^l$.

When $\zz^\nu = z_1^{\nu_1} \cdots z_l^{\nu_l}$ for a dominant weight
$\nu$ (a weight $\nu \in \ZZ^l$ is dominant if $\nu_1 \ge \cdots \ge
\nu_l$), $\sigmabold (\zz ^{\nu }) = \chi _{\nu }$ is an irreducible
$\GL _{l}$ character.  For an arbitrary weight $\gamma \in \ZZ ^{l}$,
either $\sigmabold (\zz ^{\gamma }) = \pm \chi _{\nu }$ for a suitable
dominant weight $\nu $, or $\sigmabold (\zz ^{\gamma }) = 0$.  We
extend $\sigmabold$ to an operator on formal $\kk$-linear combinations
$\sum_{\gamma \in \ZZ^l} c_\gamma \zz^\gamma$ by applying it term by
term, giving an infinite formal linear combination of irreducible
$\GL_l$ characters $\sum_{\nu} a_\nu \chi_\nu = \sum_{\gamma \in \ZZ^l}
c_\gamma \sigmabold(\zz^\gamma) $.
This makes sense because for each dominant weight $\nu $, the
set of monomials $\zz ^{\gamma }$ such that $\sigmabold (\zz ^{\gamma
}) = \pm \chi _{\nu }$ is finite.

Recall that the {\em polynomial characters} of $\GL _l$ are the
irreducible characters $\chi _{\nu }$ for which $\nu $ is a partition,
that is, $\nu _l\geq 0$.  Given any formal $\kk$-linear combination
$\sum _{\nu} a_{\nu } \chi _{\nu }$ of irreducible $\GL _{l}$
characters, we define its \emph{polynomial truncation} by
\begin{align}
\label{ed pol}
\pol_X \! \Big(\sum _{\nu}a_{\nu }\spa \chi _{\nu }\Big)
= \sum _{ \nu_l \ge 0 }a_{\nu }\spa s_{\nu }(X).
\end{align}
In principle the right hand side is an infinite formal sum of
symmetric functions, but, for instance, if $\sum a_{\nu } \chi _{\nu
}$ is homogeneous of degree $d$, then the right hand side is an
ordinary symmetric function, homogeneous of degree $d$.

We define a related operator $\mathbf{h}_X$
on Laurent polynomials $f(\zz )$
by
\begin{equation}\label{e:Weyl-symmetrizer 2}
\mathbf{h}_X(f(\zz)) =
\pol_X  \sigmabold
\bigg(\frac{f(\zz)}{\prod_{\alpha_{ij} \in R_+} (1-z_i/z_j)} \bigg),
\end{equation}
where the factors $(1- z_{i}/z_{j})^{-1}$ are expanded as
geometric series in $z_i/z_j$ before applying  $\sigmabold$.
When $f$ is a monomial, it is well known \cite{Thomas81} that
\begin{align}
\label{e hX equals h}
\mathbf{h}_X(\zz^\gamma) &= h_\gamma(X),
\end{align}
where for any integer vector $\gamma \in \ZZ^l$, we define $h_\gamma =
h_{\gamma_1}\cdots h_{\gamma_l}$ to be the product of complete
homogeneous symmetric functions, with $h_d$ for $d\leq 0$ interpreted
as $h_0 =1$, or $h_d = 0$ for $d< 0$.

We again extend the definition to formal linear combinations of
monomials, so that $\mathbf{h}_X(\sum_{\gamma \in \ZZ^l} c_\gamma
\zz^\gamma) = \sum_{\gamma \in \ZZ^l} c_\gamma h_\gamma (X)$.
With this interpretation, \eqref{e:Weyl-symmetrizer 2} still remains
valid when $f$ is a power series in $z_i/z_j$ for $i<j$.
As with $\pol_X$, in principle $\mathbf{h}_X(\sum_{\gamma} c_\gamma
\zz^\gamma)$ is an infinite formal sum of symmetric functions, but,
for instance, if $\sum_{\gamma} c_\gamma \zz^\gamma $ is homogeneous
of degree $d$, then $\mathbf{h}_X(\sum_{\gamma } c_\gamma \zz^\gamma)$
is an ordinary symmetric function, homogeneous of degree $d$.

\begin{remark}\label{rem:geometric-series}
Below we will write other formulas involving $\sigmabold $ applied to
an expression with denominator factors resembling those in
\eqref{e:Weyl-symmetrizer 2}.  Our convention is always to expand
denominator factors of the form $(1-c \spa z_i/z_j)$ for $c \in \kk$
and $i < j$ as geometric series $(1 - c z_i/z_j)^{-1} = 1 + c z_i/z_j
+ \cdots$ before applying $\sigmabold $.
\end{remark}

\subsection{Raising operator formulas for modified Hall-Littlewood
polynomials} \label{ss:HL-raising}

To set the stage for our raising operator formula for modified
Macdonald polynomials, we review two different raising operator
formulas for the modified Hall-Littlewood polynomials.  Both formulas
naturally reflect the geometry of the flag variety $G/B$; one realizes
$\Htild_\mu(X;0,t)$ as the graded Euler characteristic of the
cotangent bundle of $G/B$ twisted by a line bundle of weight $-\mu$,
while the other is the graded Euler characteristic of the cotangent
bundle of $G/P_\mu$, where $P_\mu$ is the parabolic subgroup whose
block sizes are the parts of $\mu $.  See \cite{BroerNormality} and
\cite{ShimWeym00} for details.

The first raising operator formula for $\Htild_\mu(X;0,t)$,
is as follows (see \cite[III
(6.3)]{Macdonald95} or \cite[(4.28) and
\S2]{Milneclassicalpartitionfunctions}):
\begin{align}
\label{e HL formula from mu} t^{\mathsf{n}(\mu)}\Htild_\mu(X;0,t^{-1}) = \pol_X
\! \Big( \sigmabold\Big(\frac{\zz^\mu}{\prod_{\alpha_{ij} \in R_+}(1 -
t \spa z_i/z_j)}\Big)\Big)\,,
\end{align}
where the denominator factors are expanded as geometric series
in accordance with Remark~\ref{rem:geometric-series}.

A second raising operator formula follows from the work of Weyman and
Shimozono-Weyman (see \cite[Theorem 6.10]{Weyman89} and \cite[\S2.3
(2) and (2.3)--(2.5)]{ShimWeym00}).  In this formula, the input
partition $\mu$ appears in the set of roots, instead of in the weight
$\zz^\mu$, as it does in formula~\eqref{e HL formula from mu}.  Given
a partition $\mu$ of $l$, consider the set partition of $\{1,\ldots,l
\}$ into intervals of lengths $\mu_{\ell(\mu)}, \dots, \mu_1$, and let
$B_\mu$ denote the set of roots $\alpha_{ij}$ such that $i < j$ appear
in distinct blocks of this partition.  Then
\begin{align}
\label{e HL formula from 1s} \Htild_\mu(X;0,t) = \omega \pol_X \!
\Big( \sigmabold\Big(\frac{z_1\cdots z_l}{\prod_{\alpha_{ij} \in
B_\mu}(1 - t \spa z_i/z_j)}\Big) \Big).
\end{align}
Here we chose to take the parts of $\mu $ in reverse order for
compatibility with the formula \eqref{et:mac formula a} given later,
but the order doesn't actually matter in \eqref{e HL formula from 1s}.

\begin{remark}\label{r:raising}
Formulas such as \eqref{e HL formula from mu} and
\eqref{e HL formula from 1s} are  traditionally written using an informal notation---as in
\cite{Macdonald95},
\cite[(4.28)]{Milneclassicalpartitionfunctions}, or \cite[\S 2]{Tamvakis11}---in which formula \eqref{e HL formula from mu}, for example,
would be written as
\begin{align}\label{et:raising op HL}
t^{\mathsf{n}(\mu)}\Htild_\mu(X;0,t^{-1}) = \frac{ 1 } { \prod_{\alpha_{i j} \in
R_+} \big(1-t \spa \mathbf{R}_{ij}\big)} \cdot s_{\mu}\,,
\end{align}
with raising
``operators'' $\mathbf{R}_{ij}$ which act on the subscript of a Schur
function $s_\gamma$ by $\mathbf{R}_{ij} \gamma = \gamma + \epsilon_i -
\epsilon_j$.  Here, Schur functions indexed by non-partition weights
are defined by $s_\gamma(X) = \pol _{X}\sigmabold(\zz^\gamma)$, which
is equal to 0 or to $\pm 1$ times a Schur function of partition
weight. Note that all the raising operators must be applied before
converting Schur functions of non-partition weights to ones indexed by
partition weights.  The $\mathbf{R}_{ij}$ are not true operators,
e.g. $\mathbf{R}_{23}s_{(1,1,1)} = \pol _{X}\sigmabold(z_1z_2^2)= 0$
but $\mathbf{R}_{12}\mathbf{R}_{23}s_{(1,1,1)} = s_{(2,1,0)}\not =0$,
so we think of \eqref{et:raising op HL} as a convenient but informal
version of \eqref{e HL formula from mu}.
\end{remark}

\section{Raising operator formula for modified Macdonald polynomials}

\subsection{The formula}
\label{ss formula for partitions}

Let $\mu$ be a partition of $l$.  The \emph{reading order} is the
total order on the boxes of $\mu$ given by $(i,j)\rdless (i',j')$ if
$j > j'$, or $j = j'$ and $i < i'$.  We let $\bx{\mu}{1}, \dots,
\bx{\mu}{l}$ denote the boxes of $\mu$ in increasing reading order,
which is the list of boxes of $\mu$ read by rows from left to right
starting from the top row, as shown in Example \ref{ex macanimal}.
For a box $b=(i,j)$, $\south(b)=(i,j-1)$ denotes the box immediately
south of $b$.  Define subsets of $R_+ = R_+(\GL_l)$ by
\begin{align}
\label{ed:Rmu}
R_{\mu} =  \big\{ \alpha_{ij} \in R_+ : \spa \south(\bx{\mu}{i}) \in \mu, \,
\south(\bx{\mu}{i}) \spa \rdlesseq \spa \bx{\mu}{j} \big\},  \\
\widehat{R}_{\mu} = \big\{ \alpha_{ij} \in R_+ : \spa \south(\bx{\mu}{i})
\in \mu, \, \south(\bx{\mu}{i}) \spa \rdless \spa \bx{\mu}{j} \big\}.
\label{ed:Rmu2}
\end{align}

\begin{defn}
\label{d mac formula} For a partition $\mu$ of $l$, define the
\emph{Macdonald series} by
\begin{align}\label{et:mac formula a}
\HS_{\mu}(\zz;q,t) = \sigmabold \Bigg(
\frac{
\prod_{\alpha_{ij} \in R_\mu \setminus \widehat{R}_\mu }
\big(1- q^{\arm(\bx{\mu}{i})+1}\spa t^{-\leg(\bx{\mu}{i})} z_i/z_j \big)
\prod_{\alpha_{ij} \in \widehat{R}_{\mu}}
\big(1-q\spa  t\, z_i/z_j \big)} {\prod_{\alpha_{ij} \in R_+} \big(1-q \, z_i/z_j\big)
\prod_{\alpha_{ij} \in R_\mu} \big(1-t \, z_i/z_j\big)}\Bigg),
\end{align}
which we interpret as an infinite formal linear combination of
irreducible $\GL _{l}$ characters by expanding the denominator factors
as geometric series, in accordance with
Remark~\ref{rem:geometric-series}.
\end{defn}

We now have the following raising operator formula for the modified
Macdonald polynomials $ \Htild_{\mu }(X; q,t)$.  The proof will be
given in \S \ref{s proof of theorems}.

\begin{thm}
\label{t mac formula}
For any partition  $\mu$ of  $l$,
\begin{align}\label{et:mac formula b}
\Htild_\mu(X;q,t) =\omega \pol_X \! \big(z_1\cdots z_l \,
\HS_\mu(\zz;q, t) \big),
\end{align}
where $\pol_X$ is as defined in \eqref{ed pol}.
\end{thm}

\begin{example}
\label{ex macanimal} For $\mu = (2,2,1,1)$, $l = 6$, the series
$\HS_{\mu}$ is visualized below with the subsets of roots $R_\mu$ and
$\widehat{R}_\mu$ drawn in an $\ell\times \ell$ grid, labeled by
matrix-style coordinates and specified by the legend, and with large
circles marking the presence of the factors involving arm and leg,
which are $(1-q \spa z_1/z_2)$, $(1-qt^{-1}z_2/z_3)$,
$(1-q^2t^{-2}z_3/z_5)$, $(1-q \spa z_4/z_6)$.
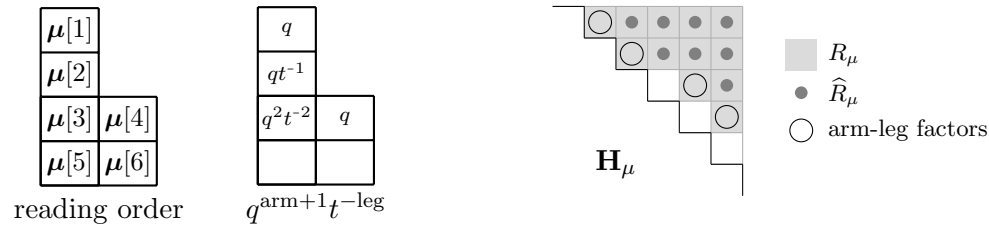
\begin{align*}
\raisebox{-5mm}{
\begin{tikzpicture}[xscale=.77, yscale = .59]
\draw[thick] (0,0) grid (1,4);
\draw[thick] (0,0) grid (2,2);
\node at (0.5, 3.5) { \footnotesize $\bx{\mu}{1}$};
\node at (0.5, 2.5) { \footnotesize $\bx{\mu}{2}$};
\node at (0.5, 1.5) { \footnotesize $\bx{\mu}{3}$};
\node at (1.5, 1.5) { \footnotesize $\bx{\mu}{4}$};
\node at (0.5, 0.5) { \footnotesize $\bx{\mu}{5}$};
\node at (1.5, 0.5) { \footnotesize $\bx{\mu}{6}$};
\node at (1,-.57) {\small reading order };
\end{tikzpicture}}
\ \ \
\raisebox{-5mm}{
\begin{tikzpicture}[xscale = .77, yscale=.59]
\draw[thick] (0,0) grid (1,4);
\draw[thick] (0,0) grid (2,2);
\node at (0.5, 3.5) {\tiny $ q $};
\node at (0.5, 2.5) {\tiny $ q t^{\shortminus 1} $};
\node at (0.5, 1.5) {\tiny $ q^2 t^{\shortminus 2} $};
\node at (1.5, 1.5) {\tiny $ q  $};
\node at (0.5, 0.5) {\tiny $ $};
\node at (1.5, 0.5) {\tiny $ $};
\node at (1,-.53) {\small $q^{\rm arm+1}t^{-{\rm leg}}$};
\end{tikzpicture}}
\quad \quad  \quad \quad \quad
\begin{tikzpicture}[scale = .42]
\draw[draw = none, fill = black!14] (2,-1) rectangle (3,-2);
 \draw[draw = none, fill = black!14] (3,-1) rectangle (4,-2);
 \draw[draw = none, fill = black!14] (3,-2) rectangle (4,-3);
 \draw[draw = none, fill = black!14] (4,-1) rectangle (5,-2);
 \draw[draw = none, fill = black!14] (4,-2) rectangle (5,-3);
 \draw[draw = none, fill = black!14] (5,-1) rectangle (6,-2);
 \draw[draw = none, fill = black!14] (5,-2) rectangle (6,-3);
 \draw[draw = none, fill = black!14] (5,-3) rectangle (6,-4);
 \draw[draw = none, fill = black!14] (6,-1) rectangle (7,-2);
 \draw[draw = none, fill = black!14] (6,-2) rectangle (7,-3);
 \draw[draw = none, fill = black!14] (6,-3) rectangle (7,-4);
 \draw[draw = none, fill = black!14] (6,-4) rectangle (7,-5);
 \draw[draw = none, fill = gray!100] (3+0.5, -1-0.5) circle (.2);
\draw[draw = none, fill = gray!100] (4+0.5, -1-0.5) circle (.2);
\draw[draw = none, fill = gray!100] (4+0.5, -2-0.5) circle (.2);
\draw[draw = none, fill = gray!100] (5+0.5, -1-0.5) circle (.2);
\draw[draw = none, fill = gray!100] (5+0.5, -2-0.5) circle (.2);
\draw[draw = none, fill = gray!100] (6+0.5, -1-0.5) circle (.2);
\draw[draw = none, fill = gray!100] (6+0.5, -2-0.5) circle (.2);
\draw[draw = none, fill = gray!100] (6+0.5, -3-0.5) circle (.2);
 \draw[draw = black!100, fill = none] (2+0.5, -1-0.5) circle (.36);
\draw[draw = black!100, fill = none] (3+0.5, -2-0.5) circle (.36);
\draw[draw = black!100, fill = none] (5+0.5, -3-0.5) circle (.36);
\draw[draw = black!100, fill = none] (6+0.5, -4-0.5) circle (.36);
\draw[thin, black!31] (1,-1) -- (7,-1);
\draw[thin, black!31] (2,-1) -- (2,-1);
\draw[thin, black!31] (2,-2) -- (7,-2);
\draw[thin, black!31] (3,-2) -- (3,-1);
\draw[thin, black!31] (3,-3) -- (7,-3);
\draw[thin, black!31] (4,-3) -- (4,-1);
\draw[thin, black!31] (4,-4) -- (7,-4);
\draw[thin, black!31] (5,-4) -- (5,-1);
\draw[thin, black!31] (5,-5) -- (7,-5);
\draw[thin, black!31] (6,-5) -- (6,-1);
\draw[thin, black!31] (6,-6) -- (7,-6);
\draw[thin, black!31] (7,-6) -- (7,-1);
\draw[thin] (1,-1) -- (2,-1);
\draw[thin] (2,-1) -- (2,-2);
\draw[thin] (2,-2) -- (3,-2);
\draw[thin] (3,-2) -- (3,-3);
\draw[thin] (3,-3) -- (4,-3);
\draw[thin] (4,-3) -- (4,-4);
\draw[thin] (4,-4) -- (5,-4);
\draw[thin] (5,-4) -- (5,-5);
\draw[thin] (5,-5) -- (6,-5);
\draw[thin] (6,-5) -- (6,-6);
\draw[thin] (6,-6) -- (7,-6);
\draw[thin] (7,-6) -- (7,-7);
\node at (6/2,-12/2) {\small $\HS_{\mu}$};
\end{tikzpicture}
\ \ \
\raisebox{6mm}{
\begin{tikzpicture}[scale = .42]
\begin{scope}
\draw[draw = none, fill = black!14] (2,-1) rectangle (3,-2);
\node[anchor = west] at (3,-1.5) {\scriptsize   $R_\mu$};
\end{scope}
\begin{scope}[yshift = -34*1]
\draw[thin, black!0] (2,-1) -- (3,-1);
\draw[thin, black!0] (2,-1) -- (2,-2);
\draw[thin, black!0] (2,-2) -- (3,-2);
\draw[thin, black!0] (3,-2) -- (3,-1);
\draw[draw = none, fill = gray!100] (2+0.5, -1-0.5) circle (.2);
\node[anchor = west] at (3,-1.5) {\scriptsize    $\widehat{R}_\mu$};
\end{scope}
\begin{scope}[yshift = -34*2]
\draw[draw = black!100, fill = none] (2+0.5, -1-0.5) circle (.36);
\node[anchor = west] at (3,-1.5) {\scriptsize    arm-leg factors};
\end{scope}
\end{tikzpicture}}
\end{align*}
\end{example}

Formula \eqref{et:mac formula b}, written in the raising operator
notation of Remark~\ref{r:raising}, becomes the formula previewed in
the introduction:
\begin{align} \label{et:raising op main}
\omega \Htild_\mu (X;q,t) = \frac{\prod_{\alpha_{ij} \in R_\mu
\setminus \widehat{R}_\mu } \big(1- q^{\arm(\bx{\mu}{i})+1}\spa
t^{-\leg(\bx{\mu}{i})} \mathbf{R}_{ij} \big) \prod_{\alpha_{ij} \in
\widehat{R}_{\mu}} \big(1-q\spa t\, \mathbf{R}_{ij} \big)}
{\prod_{\alpha_{ij} \in R_+} \big(1-q \, \mathbf{R}_{ij}\big)
\prod_{\alpha_{ij} \in R_\mu} \big(1-t \, \mathbf{R}_{ij}\big)} \cdot
s_{1^l}\,,
\end{align}
where $1^l$ denotes the all 1's vector of length $l$.

\subsection{Specializations} \label{ss specialization}

It is instructive to see how the well-known specializations of
Macdonald polynomials, $\Htild _{\mu }(X;1,1)=e_{1}^{|\mu |}(X)$ and
the Hall-Littlewood specialization $\Htild _{\mu }(X;0,t)$, can be
recovered from formula~\eqref{et:mac formula b}.

First we consider the specialization $q=t=1$.  After specializing, the
arm-leg and $(1-q\,t \,z_i/z_j)$ numerator factors in the definition
\eqref{et:mac formula a} of $\HS_\mu$ cancel with the $(1-t
\,z_i/z_j)$ factors in the denominator.  Hence
\begin{align}\label{e:11}
\HS_\mu(\zz;1,1) = \sigmabold\Big(\frac{1}{\prod_{\alpha_{ij} \in
R_+}(1 - z_i/z_j)}\Big).
\end{align}
Then, using \eqref{e:Weyl-symmetrizer 2} and \eqref{e hX equals h}, we
recover the specialization
\begin{equation}\label{e:q=t=1-specialization}
\Htild_\mu(X;1,1) = \omega \pol
_{X}(z_1\cdots z_l \spa \HS_\mu(\zz;1,1)) = e_{1^l} = e_1^l.
\end{equation}

Now consider the specialization $q=0$.  After specializing
\eqref{et:mac formula a}, all numerator factors and the $(1-q
z_i/z_j)$ denominator factors reduce to 1.  This gives
\begin{align}\label{e q0}
z_1\cdots z_l \, \HS_\mu(\zz;0,t) = \sigmabold\Big(\frac{z_1\cdots
z_l}{\prod_{\alpha_{ij} \in R_\mu}(1 - t \spa z_i/z_j)}\Big).
\end{align}
Let $B_\mu$ be the set of roots defined before \eqref{e HL formula
from 1s}, which can also be described as the set of positive roots
above a block diagonal matrix with block sizes $\mu_{\ell(\mu)},
\dots, \mu_1$.  Now $R_\mu$ is contained in $B_\mu$ and $B_\mu
\setminus R_\mu$ consists of triangular regions of roots between each
pair of consecutive blocks.  We reach the Hall-Littlewood raising
operator formula using the identity
\begin{align}\label{e triangular to block}
\sigmabold\Big(\frac{z_1\cdots z_l}{\prod_{\alpha_{ij} \in B_\mu}(1 -
t \spa z_i/z_j)}\Big)
 = \sigmabold\Big(\frac{z_1\cdots z_l}{\prod_{\alpha_{ij} \in R_\mu}(1
- t \spa z_i/z_j)}\Big).
\end{align}
This is proven by removing these triangular regions from $B_\mu$ one
root at a time (starting with the bottommost region), and using the
following simplified version of
\cite[Lemma 8.9]{BMPS2} to show that the corresponding functions
remain the same at each step.

\begin{lemma}
\label{l root expansion} Let $k \in \{1,\dots, l-1\}$ and suppose that
$B \subseteq R_+(\GL_l)$ is a set of roots such that
$\prod_{\alpha_{ij}\in B} (1-tz_i/z_j)$ is fixed by the simple
reflection $s_k$.  If $B$ contains a root $\alpha = \alpha_{k+1,j}$ for some $j
> k+1$, then
\begin{align}\label{e root expansion}
\sigmabold\Big(\frac{z_1\cdots z_l}{\prod_{\alpha_{ij} \in B}(1 - t
\spa z_i/z_j)}\Big)
 = \sigmabold\Big(\frac{z_1\cdots z_l}{\prod_{\alpha_{ij} \in B
\setminus \alpha}(1 - t \spa z_i/z_j)}\Big).
\end{align}
\end{lemma}

Combining \eqref{e q0} and \eqref{e triangular to block} now shows
that the raising operator formula
\eqref{et:mac formula b} for modified Macdonald polynomials reduces at
$q=0$ to the raising operator formula
\eqref{e HL formula from 1s} for Hall-Littlewood polynomials.

\subsection{A formula for any column diagram}
\label{ss:column-diagram-formula}

We generalize Theorem \ref{t mac formula} to raising operator formulas
for $ \Htild_\mu(X;q,t)$ in terms of any diagram obtained by
rearranging the columns of $\mu$.  For $\beta = (\beta_1, \dots,
\beta_k) \in \ZZ_{+}^k$, we define the \emph{column diagram} of
$\beta$ to be the set
\begin{align}\label{e diagram beta}
\dgrm{\beta} = \big\{(i,j) \in \ZZ_{+}^2 : \spa 1 \le i \le k, \, 1 \le j
\le \beta_i \big\}.
\end{align}
In particular, if $\beta $ is the transpose of a partition $\mu $, then $\beta^*$ is the diagram of $\mu $.

All of the definitions associated to $\mu$ carry over to
$\dgrm{\beta}$ essentially unchanged, except for the $\arm$ statistic.
As before, we identify $(i,j)$ with the lattice square or box whose
northeast corner has coordinates $(x,y) = (i,j)$; we again say that
this box is in column $i$ and row $j$.  Define the reading order
$\rdless$ on $\dgrm{\beta}$ by $(i,j) \rdless (i',j')$ if $j > j'$, or
$j = j'$ and $i < i'$; for the boxes $\bx{\beta}{1},\dots,
\bx{\beta}{l}$ of $\dgrm{\beta}$ listed in increasing reading order,
we define subsets of  $R_+ = R_+(\GL_l)$ by
\begin{align}\label{ed:Rbeta}
R_{\beta^*} = \big\{ \alpha_{ij} \in R_+  :  \spa \south(\bx{\beta}{i}) \in
\dgrm{\beta}, \,
\south(\bx{\beta}{i}) \spa \rdlesseq \spa \bx{\beta}{j} \big\},  \\
\widehat{R}_{\beta^*} = \big\{ \alpha_{ij} \in R_+  : \spa
\south(\bx{\beta}{i}) \in \dgrm{\beta}, \, \south(\bx{\beta}{i}) \spa \rdless
\spa \bx{\beta}{j} \big\}.  \label{ed:Rbeta2}
\end{align}
For $\beta \in \ZZ_{+}^k$ and a box $b = (i,j) \in \dgrm{\beta}$ with  $j > 1$,
define the \emph{arm} and \emph{leg} of $b$ by
\begin{align}
\label{ed leg}
\leg(b)  &= \beta_i -j  = \text{(number of boxes strictly  north of  $b$)},\\
\label{ed arm}
\arm(b) &= \big| \big\{i' \in \{1,\dots, i-1\} : j\!-\!1 \le \beta_{i'} <
\beta_i \big\} \sqcup \!  \big\{i' \in \{i+1,\dots, k\} : j \le
\beta_{i'} \le \beta_i \big\} \big|.\!\!\!\!\!
\end{align}
In words, $\arm (b)$ is the number of boxes strictly
east of $b$ in columns of height $\beta _{i'}\leq \beta _{i}$ or
strictly west of $\south (b)$ in columns of height $\beta _{i'}<\beta_{i}$.

\begin{example}
For $\beta = (3,2,4,3,4,2,1,3)$, the column diagram of $\beta$ is
displayed below, along with the $\arm$ of box $\bullet=(4,2)$
where
the $a$'s mark the boxes contributing to $\arm(\bullet)$.
\begin{align}\label{e:column-diagram}
\dgrm{\beta} \spa =
\raisebox{-7mm}{
\begin{tikzpicture}[xscale = .38, yscale = .38]
\draw[thick] (0,0) grid (1,3);
\draw[thick] (1,0) grid (2,2);
\draw[thick] (2,0) grid (3,4);
\draw[thick] (3,0) grid (4,3);
\draw[thick] (4,0) grid (5,4);
\draw[thick] (5,0) grid (6,2);
\draw[thick] (6,0) grid (7,1);
\draw[thick] (7,0) grid (8,3);
\node at (1.5, 0.5) {$a$};
\node at (5.5, 1.5) {$a$};
\node at (7.5, 1.5) {$a$};
\node at (3.5, 1.46) { $\bullet$};
\end{tikzpicture}
}
\quad  \quad  \ \arm(\bullet) = 3
\end{align}
\end{example}

\begin{thm}\label{t mac formula 2}
Given $\beta \in \ZZ_{+}^k$, using the notation just above, define the
\emph{Macdonald series}
\begin{align}\label{et:mac formula a 2}
\HS_{\beta^*}(\zz;q,t) = \sigmabold \Bigg( \frac{ \prod_{\alpha_{ij} \in
R_{\beta^*} \setminus \widehat{R}_{\beta^*}} \big(1-
q^{\arm(\bx{\beta}{i})+1}\spa t^{-\leg(\bx{\beta}{i})} z_i/z_j \big)
\prod_{\alpha_{ij} \in \widehat{R}_{\beta^*}} \big(1-q\spa t\, z_i/z_j
\big)} {\prod_{\alpha_{ij} \in R_+} \big(1-q \, z_i/z_j\big)
\prod_{\alpha_{ij} \in R_{\beta^*}} \big(1-t \, z_i/z_j\big)}\Bigg),
\end{align}
which we regard as an infinite formal linear combination of
irreducible $\GL _{l}$ characters using the convention of Remark
\ref{rem:geometric-series}.  Then, for any partition $\mu $ and any
rearrangement $\beta $ of $\mu ^*$, we have
\begin{equation}
    \Hbold_{\beta^*}(\zz;q,t) =  \Hbold_\mu(\zz;q,t).
\end{equation}
\end{thm}

The proof will be given in \S \ref{s:m,n-Macdonalds}.  Combining
Theorems~\ref{t mac formula} and \ref{t mac formula 2} gives the following
generalization of \eqref{et:mac formula b}.

\begin{cor}
\label{cor:mac formula b 2}
For any partition  $\mu$ of  $l$
and any rearrangement $\beta $ of $\mu ^*$,
\begin{align}\label{et:mac formula b 2}
\Htild_\mu(X;q,t) = \omega \pol _{X}\! \big(z_1\cdots z_l \,
\HS_{\beta^*}(\zz;q, t) \big).
\end{align}
\end{cor}

\begin{remark}\label{rem:mac-2-implies-mac}
Theorem \ref{t mac formula} is the special case of
Corollary~\ref{cor:mac formula b 2} when $\beta$ is a partition, since
the formula \eqref{et:mac formula a} for
$\Hbold_\mu(\zz;q,t)$ is exactly \eqref{et:mac formula a 2} with $\beta
= \mu^*$.
\end{remark}

\begin{example}
\label{ex macanimal2} For $\beta = (1,4,2,4)$, $l = 11$, the series
$\HS_{\beta^*}$ is visualized below in the style of Example \ref{ex
macanimal}; the factors involving arm and leg, marked by the large
circles, are $(1-q^2 \spa z_1/z_3)$, $(1-qz_2/z_4)$,
$(1-q^2t^{-1}z_3/z_5)$, $(1-q^2t^{-1} \spa z_4/z_7)$, $(1-q^4t^{-2}
\spa z_5/z_9)$, $(1-q^2 \spa z_6/z_{10})$, $(1-q^3t^{-2} \spa
z_7/z_{11})$.
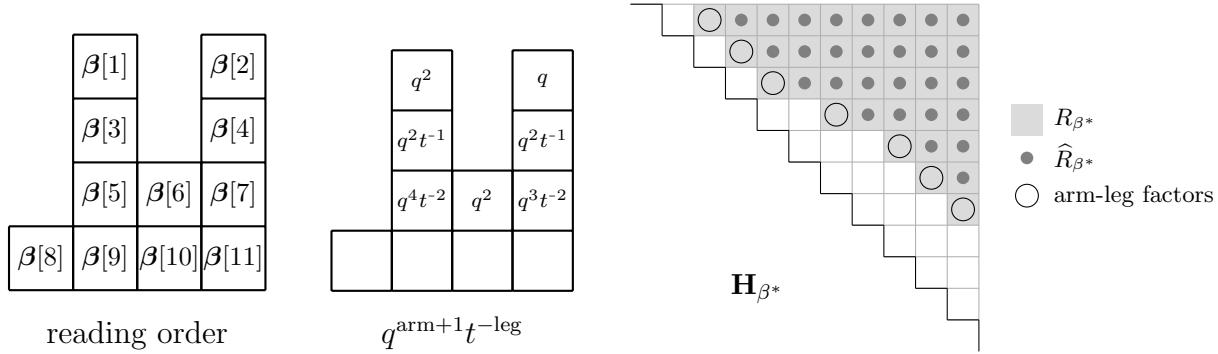
\begin{align*}
  \
  \begin{tikzpicture}[scale=0.85]
    \drawDg{1,4,2,4}
    \setcounter{boxnum}{1};
    \foreach \x\y in {2/4,4/4,2/3,4/3,2/2,3/2,4/2,1/1,2/1,3/1,4/1} {
      \node at (\x-0.5,\y-0.5) {\footnotesize $\bx{\beta}{\theboxnum}$};
      \addtocounter{boxnum}{1};
     };
     \node at (2,-.7) {reading order};
  \end{tikzpicture}
  \qquad \!
  \begin{tikzpicture}[scale=0.8]
    \drawDg{1,4,2,4}
    \node at (1.5,3.5) {\tiny $q^2$};
    \node at (1.5,2.5) {\tiny $q^2t^{\shortminus 1}$};
    \node at (1.5,1.5) {\tiny $q^4t^{\shortminus 2}$};
    \node at (2.5,1.5) {\tiny $q^2$};
    \node at (3.5,3.5) {\tiny $q$};
    \node at (3.5,2.5) {\tiny $q^2t^{\shortminus 1}$};
    \node at (3.5,1.5) {\tiny $q^3t^{\shortminus 2}$};
     \node at (2,-.7) {\(q^{\arm+1}t^{-\leg}\)};
  \end{tikzpicture}
\quad  \
\raisebox{1mm}{
\begin{tikzpicture}[scale = .42]
\draw[draw = none, fill = black!14] (3,-1) rectangle (4,-2);
 \draw[draw = none, fill = black!14] (4,-1) rectangle (5,-2);
 \draw[draw = none, fill = black!14] (4,-2) rectangle (5,-3);
 \draw[draw = none, fill = black!14] (5,-1) rectangle (6,-2);
 \draw[draw = none, fill = black!14] (5,-2) rectangle (6,-3);
 \draw[draw = none, fill = black!14] (5,-3) rectangle (6,-4);
 \draw[draw = none, fill = black!14] (6,-1) rectangle (7,-2);
 \draw[draw = none, fill = black!14] (6,-2) rectangle (7,-3);
 \draw[draw = none, fill = black!14] (6,-3) rectangle (7,-4);
 \draw[draw = none, fill = black!14] (7,-1) rectangle (8,-2);
 \draw[draw = none, fill = black!14] (7,-2) rectangle (8,-3);
 \draw[draw = none, fill = black!14] (7,-3) rectangle (8,-4);
 \draw[draw = none, fill = black!14] (7,-4) rectangle (8,-5);
 \draw[draw = none, fill = black!14] (8,-1) rectangle (9,-2);
 \draw[draw = none, fill = black!14] (8,-2) rectangle (9,-3);
 \draw[draw = none, fill = black!14] (8,-3) rectangle (9,-4);
 \draw[draw = none, fill = black!14] (8,-4) rectangle (9,-5);
 \draw[draw = none, fill = black!14] (9,-1) rectangle (10,-2);
 \draw[draw = none, fill = black!14] (9,-2) rectangle (10,-3);
 \draw[draw = none, fill = black!14] (9,-3) rectangle (10,-4);
 \draw[draw = none, fill = black!14] (9,-4) rectangle (10,-5);
 \draw[draw = none, fill = black!14] (9,-5) rectangle (10,-6);
 \draw[draw = none, fill = black!14] (10,-1) rectangle (11,-2);
 \draw[draw = none, fill = black!14] (10,-2) rectangle (11,-3);
 \draw[draw = none, fill = black!14] (10,-3) rectangle (11,-4);
 \draw[draw = none, fill = black!14] (10,-4) rectangle (11,-5);
 \draw[draw = none, fill = black!14] (10,-5) rectangle (11,-6);
 \draw[draw = none, fill = black!14] (10,-6) rectangle (11,-7);
 \draw[draw = none, fill = black!14] (11,-1) rectangle (12,-2);
 \draw[draw = none, fill = black!14] (11,-2) rectangle (12,-3);
 \draw[draw = none, fill = black!14] (11,-3) rectangle (12,-4);
 \draw[draw = none, fill = black!14] (11,-4) rectangle (12,-5);
 \draw[draw = none, fill = black!14] (11,-5) rectangle (12,-6);
 \draw[draw = none, fill = black!14] (11,-6) rectangle (12,-7);
 \draw[draw = none, fill = black!14] (11,-7) rectangle (12,-8);
 \draw[draw = none, fill = gray!100] (4+0.5, -1-0.5) circle (.2);
\draw[draw = none, fill = gray!100] (5+0.5, -1-0.5) circle (.2);
\draw[draw = none, fill = gray!100] (5+0.5, -2-0.5) circle (.2);
\draw[draw = none, fill = gray!100] (6+0.5, -1-0.5) circle (.2);
\draw[draw = none, fill = gray!100] (6+0.5, -2-0.5) circle (.2);
\draw[draw = none, fill = gray!100] (6+0.5, -3-0.5) circle (.2);
\draw[draw = none, fill = gray!100] (7+0.5, -1-0.5) circle (.2);
\draw[draw = none, fill = gray!100] (7+0.5, -2-0.5) circle (.2);
\draw[draw = none, fill = gray!100] (7+0.5, -3-0.5) circle (.2);
\draw[draw = none, fill = gray!100] (8+0.5, -1-0.5) circle (.2);
\draw[draw = none, fill = gray!100] (8+0.5, -2-0.5) circle (.2);
\draw[draw = none, fill = gray!100] (8+0.5, -3-0.5) circle (.2);
\draw[draw = none, fill = gray!100] (8+0.5, -4-0.5) circle (.2);
\draw[draw = none, fill = gray!100] (9+0.5, -1-0.5) circle (.2);
\draw[draw = none, fill = gray!100] (9+0.5, -2-0.5) circle (.2);
\draw[draw = none, fill = gray!100] (9+0.5, -3-0.5) circle (.2);
\draw[draw = none, fill = gray!100] (9+0.5, -4-0.5) circle (.2);
\draw[draw = none, fill = gray!100] (10+0.5, -1-0.5) circle (.2);
\draw[draw = none, fill = gray!100] (10+0.5, -2-0.5) circle (.2);
\draw[draw = none, fill = gray!100] (10+0.5, -3-0.5) circle (.2);
\draw[draw = none, fill = gray!100] (10+0.5, -4-0.5) circle (.2);
\draw[draw = none, fill = gray!100] (10+0.5, -5-0.5) circle (.2);
\draw[draw = none, fill = gray!100] (11+0.5, -1-0.5) circle (.2);
\draw[draw = none, fill = gray!100] (11+0.5, -2-0.5) circle (.2);
\draw[draw = none, fill = gray!100] (11+0.5, -3-0.5) circle (.2);
\draw[draw = none, fill = gray!100] (11+0.5, -4-0.5) circle (.2);
\draw[draw = none, fill = gray!100] (11+0.5, -5-0.5) circle (.2);
\draw[draw = none, fill = gray!100] (11+0.5, -6-0.5) circle (.2);
\draw[draw = black!100, fill = none] (3+0.5, -1-0.5) circle (.36);
\draw[draw = black!100, fill = none] (4+0.5, -2-0.5) circle (.36);
\draw[draw = black!100, fill = none] (5+0.5, -3-0.5) circle (.36);
\draw[draw = black!100, fill = none] (7+0.5, -4-0.5) circle (.36);
\draw[draw = black!100, fill = none] (9+0.5, -5-0.5) circle (.36);
\draw[draw = black!100, fill = none] (10+0.5,-6-0.5) circle (.36);
\draw[draw = black!100, fill = none] (11+0.5, -7-0.5) circle (.36);
\draw[thin, black!31] (1,-1) -- (12,-1);
\draw[thin, black!31] (2,-1) -- (2,-1);
\draw[thin, black!31] (2,-2) -- (12,-2);
\draw[thin, black!31] (3,-2) -- (3,-1);
\draw[thin, black!31] (3,-3) -- (12,-3);
\draw[thin, black!31] (4,-3) -- (4,-1);
\draw[thin, black!31] (4,-4) -- (12,-4);
\draw[thin, black!31] (5,-4) -- (5,-1);
\draw[thin, black!31] (5,-5) -- (12,-5);
\draw[thin, black!31] (6,-5) -- (6,-1);
\draw[thin, black!31] (6,-6) -- (12,-6);
\draw[thin, black!31] (7,-6) -- (7,-1);
\draw[thin, black!31] (7,-7) -- (12,-7);
\draw[thin, black!31] (8,-7) -- (8,-1);
\draw[thin, black!31] (8,-8) -- (12,-8);
\draw[thin, black!31] (9,-8) -- (9,-1);
\draw[thin, black!31] (9,-9) -- (12,-9);
\draw[thin, black!31] (10,-9) -- (10,-1);
\draw[thin, black!31] (10,-10) -- (12,-10);
\draw[thin, black!31] (11,-10) -- (11,-1);
\draw[thin, black!31] (11,-11) -- (12,-11);
\draw[thin, black!31] (12,-11) -- (12,-1);
\draw[thin] (1,-1) -- (2,-1);
\draw[thin] (2,-1) -- (2,-2);
\draw[thin] (2,-2) -- (3,-2);
\draw[thin] (3,-2) -- (3,-3);
\draw[thin] (3,-3) -- (4,-3);
\draw[thin] (4,-3) -- (4,-4);
\draw[thin] (4,-4) -- (5,-4);
\draw[thin] (5,-4) -- (5,-5);
\draw[thin] (5,-5) -- (6,-5);
\draw[thin] (6,-5) -- (6,-6);
\draw[thin] (6,-6) -- (7,-6);
\draw[thin] (7,-6) -- (7,-7);
\draw[thin] (7,-7) -- (8,-7);
\draw[thin] (8,-7) -- (8,-8);
\draw[thin] (8,-8) -- (9,-8);
\draw[thin] (9,-8) -- (9,-9);
\draw[thin] (9,-9) -- (10,-9);
\draw[thin] (10,-9) -- (10,-10);
\draw[thin] (10,-10) -- (11,-10);
\draw[thin] (11,-10) -- (11,-11);
\draw[thin] (11,-11) -- (12,-11);
\draw[thin] (12,-11) -- (12,-12);
\node at (10/2,-20/2) {\small $\HS_{\beta^*}$};
\end{tikzpicture}}
\ \
\raisebox{19mm}{
\begin{tikzpicture}[scale = .42]
\begin{scope}
\draw[draw = none, fill = black!14] (2,-1) rectangle (3,-2);
\node[anchor = west] at (3,-1.5) {\scriptsize   $R_{\beta^*}$};
\end{scope}
\begin{scope}[yshift = -34*1]
\draw[thin, black!0] (2,-1) -- (3,-1);
\draw[thin, black!0] (2,-1) -- (2,-2);
\draw[thin, black!0] (2,-2) -- (3,-2);
\draw[thin, black!0] (3,-2) -- (3,-1);
\draw[draw = none, fill = gray!100] (2+0.5, -1-0.5) circle (.2);
\node[anchor = west] at (3,-1.5) {\scriptsize    $\widehat{R}_{\beta^*}$};
\end{scope}
\begin{scope}[yshift = -34*2]
\draw[draw = black!100, fill = none] (2+0.5, -1-0.5) circle (.36);
\node[anchor = west] at (3,-1.5) {\scriptsize    arm-leg factors};
\end{scope}
\end{tikzpicture}
  }
\end{align*}
\end{example}

\subsection{A strengthening of Macdonald positivity}
\label{ss:positivity}

By Theorem~\ref{t mac formula}, the truncation of the series $\Hbold
_{\mu } = \Hbold _{\mu }(\zz ;q,t)$ to polynomial characters yields
the Macdonald polynomial $\Htild _{\mu }(X;q,t) =\omega \pol
_{X}(z_{1}\cdots z_{l}\, \Hbold _{\mu })$, which is known
\cite{Haiman01} to have positive coefficients $\Ktild _{\lambda ,\mu
}(q,t)\in \NN [q,t]$ in terms of Schur functions.

The full series $\Hbold _{\mu }$ is encapsulated by generalized
Macdonald polynomials $\Htild _{\mu }^{1,n}(X;q,t)$ which we define in
\S \ref{s:m,n-Macdonalds} using the action of the elliptic Hall
algebra $\Ecal $ of Burban and Schiffmann \cite{BurbSchi12} on
symmetric functions.  Based on extensive computations, we were led to
the following positivity conjecture for $\Hbold _{\mu }$, which
generalizes the positivity theorem for Macdonald polynomials.

\begin{conj}\label{conj:Hmu-positivity}
For every partition $\mu $ of $l$, the series $\Hbold_\mu(\zz;q,t)$ is
a positive sum of irreducible $\GL_l$ characters; that is, the
coefficients in
\begin{equation}\label{e:Hmu-positivity}
\Hbold_\mu(\zz;q,t) = \sum_{\nu } \Kbold _{\nu , \mu}(q,t)\, \chi_\nu
\end{equation}
are polynomials $\Kbold _{\nu , \mu}(q,t)\in \NN[q,t]$ with
non-negative integer coefficients.
\end{conj}

Example \ref{ex 1nmacs} gives some instances of these coefficients, in the notation of the equivalent Conjecture~\ref{cj 1n Schur positive}.

\section{Proof of Theorem \ref{t mac formula} and Corollary \ref{cor:mac formula b 2}}
\label{s proof of theorems}

As noted in Remark~\ref{rem:mac-2-implies-mac}, Theorem~\ref{t mac formula} is a special case of Corollary \ref{cor:mac formula b 2}.
We prove
the latter using two main ingredients: the Haglund-Haiman-Loehr
formulas \cite{HaHaLo05, HaHaLo08} and our recent formula for $\nabla$
on an LLT polynomial \cite{BHMPS-llt}.  We explain these two
ingredients after a recap of LLT polynomials.

\subsection{LLT polynomials}
\label{ss:LLT}

We recall the definition and basic properties of LLT polynomials
\cite{LaLeTh97}, using the `attacking inversions' formulation from
\cite{HaHaLoReUl05}.

A {\it skew diagram} is a difference $\nu = \lambda / \mu$
of partition diagrams $\mu\subseteq\lambda$.
The {\em content} of a box $b = (i,j)$ in row $j$, column $i$ of a
(skew) diagram is $c(b) = i-j$.

Let $\nubold = (\nu_{(1)},\ldots,\nu_{(k)})$ be a tuple of skew
diagrams.  We consider the set of boxes in $\nubold $ to be the
disjoint union of the sets of boxes in the $\nu _{(i)}$, and define
the {\em adjusted content} of a box $a\in \nu_{(i)}$ to be $\ctild (a)
= c(a)+i\epsilon $, where $\epsilon $ is a fixed positive number such
that $k \epsilon <1$.

A {\em diagonal} in $\nubold $ is the set of boxes of a fixed adjusted
content, that is, a diagonal of fixed content in one of the $\nu_{(i)}$.

The {\em reading order} on $\nubold $ is the total ordering $<$ on the
boxes of $\nubold $ such that $a<b \Rightarrow \ctild (a)\leq
\ctild(b)$ and boxes on each diagonal increase from southwest to
northeast.
See Example \ref{ex:HHL full example}.  An {\em attacking
pair} is an ordered pair of boxes $(a,b)$ in $ \nubold$ such that
$a<b$ in reading order and $0<\ctild(b)-\ctild(a)<1$.

A {\em semistandard tableau} on the tuple $\nubold $ is a map $T\colon
\nubold \rightarrow \ZZ _{+}$ which restricts to a semistandard Young
tableau on each component $\nu_{(i)}$.  The set of these is denoted
$\SSYT (\nubold )$.  An {\em attacking inversion} in $T$ is an
attacking pair $(a,b)$ such that $T(a)>T(b)$.  The number of attacking
inversions in $T$ is denoted $\inv (T)$.

\begin{defn}\label{def:G-nu}
The {\em LLT polynomial} indexed by a tuple of skew diagrams $\nubold
$ is the generating function
\begin{equation}\label{e:G-nu}
\Gcal_{\nubold  }(X;q) = \sum _{T\in \SSYT (\nubold)}q^{\inv (T)}\xx ^{T},
\end{equation}
where $\xx ^{T} = \prod _{a\in \nubold  } x_{T(a)}$.
By
\cite{HaHaLoReUl05,LaLeTh97},  $\Gcal_{\nubold  }(X;q) $  is known to be symmetric.
\end{defn}

\subsection{The Haglund-Haiman-Loehr formula}\label{ss HHL formula}
\label{ss:HHL}

Haglund-Haiman-Loehr \cite{HaHaLo05} gave a formula for the modified
Macdonald polynomials $H_\mu(X;q,t)$ as a positive sum of LLT
polynomials, and this was generalized in \cite{HaHaLo08} to give many
different expressions for $H_\mu(X;q,t)$ as a positive sum of LLT
polynomials, one for each rearrangement $\beta$ of $\mu^*$.  We now
recall this formula.

A \emph{ribbon} is a connected skew shape containing no $2 \times 2$
block of boxes.

For $\beta \in \ZZ_{+}^k$, let
$V_\beta = \{(\bx{\beta}{i},\bx{\beta}{j})  :  \bx{\beta}{j} =
\south(\bx{\beta}{i}) \}$
be the set of ordered pairs of boxes that form vertical
dominoes in $\dgrm{\beta}$.

\begin{defn}
\label{d:nuofS}
For each $S
\subseteq V_\beta$, define
$\nubold(S) = (\nu_{(1)},\dots,\nu_{(k)})$
to be the $k$-tuple of ribbons where the $i$-th ribbon $\nu_{(i)}$ is
determined by
\begin{itemize}
  \setlength\itemsep{1ex}
    \item [(i)]  $\nu_{(i)}$ has $\beta_i$ boxes, of contents
      $-1,-2,\dots, -\beta_i$, and
    \item [(ii)]  the boxes of contents $-j$ and $-j+1$ in
       $\nu_{(i)}$ form a vertical domino  if and only if
    the domino $((i,j),(i,j-1))$ in $V_\beta$ belongs to $S$.
\end{itemize}
\end{defn}

\begin{thm}[{\cite{HaHaLo05,HaHaLo08}}] \label{t HHL formula}
Let $\mu$ be a partition and let $\beta$ be any rearrangement of
$\mu ^{*}$.
 Then
\begin{align}
\label{eq:HHL formula} \Htild_\mu(X;q,t) = \sum_{S \subseteq V_\beta}
\Big( \prod_{(\bx{\beta}{i},\bx{\beta}{j}) \in S}
q^{-\arm(\bx{\beta}{i})} \spa t^{\leg(\bx{\beta}{i})+1} \Big)
\Gcal_{\nubold(S)}(X;q).
\end{align}
\end{thm}

Theorem \ref{t HHL formula} in the case that $\beta$ is a partition is immediate from Theorem 2.2, Equation 23, and Proposition 3.4 in
\cite{HaHaLo05}, while the generalization to any composition $\beta$
is addressed in \cite[Theorem 5.1.1]{HaHaLo08}.  Note that our
conventions for diagrams, arms, and legs are the same as those in
\cite{HaHaLo08} except that we have reversed the order of the columns
(which makes our notation consistent with that in \cite{HaHaLo05});
we have also used the symmetry property $H_\mu(X;q,t) =
H_{\mu^*}(X;t,q)$ to translate from the exact version stated in
\cite[Theorem 5.1.1]{HaHaLo08}.

\begin{example}
\label{ex:HHL full example} For $\beta = (2,3)$, the 8 summands
appearing on the right hand side of \eqref{eq:HHL formula} are
illustrated by drawing $\nubold(S)$ with boxes labeled in reading
order and with the corresponding coefficient
$\prod_{(\bx{\beta}{i},\bx{\beta}{j}) \in S} q^{-\arm(\bx{\beta}{i})}
\spa t^{\leg(\bx{\beta}{i})+1}$ beside it; the vertical dominoes of
\(\dgrm{\beta}\) are denoted \(v_1 = (\bx{\beta}{1},\bx{\beta}{3}),
v_2=(\bx{\beta}{2},\bx{\beta}{4}),\) \(v_3 =
(\bx{\beta}{3},\bx{\beta}{5})\).  The arm and leg statistics for
$\dgrm{\beta}$ are shown on the left.

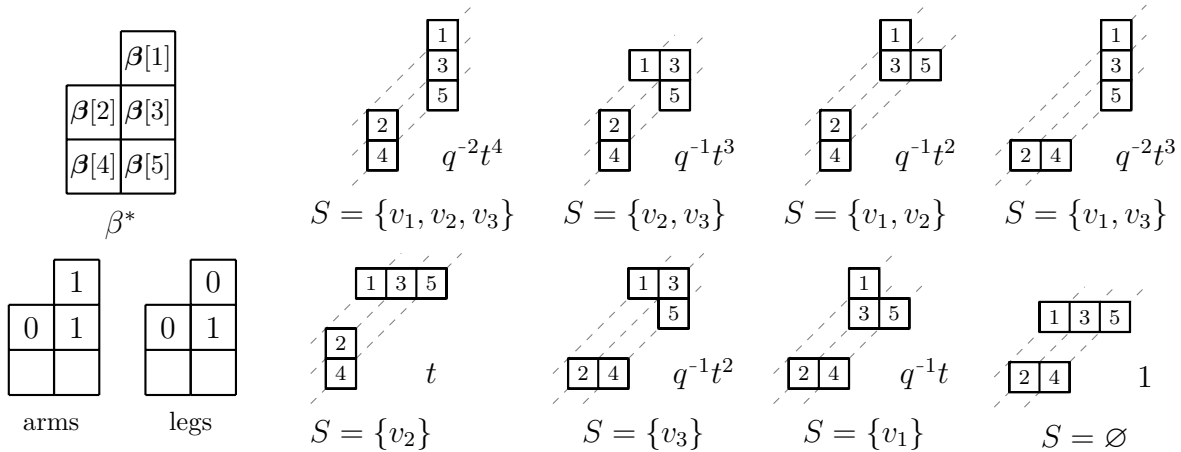
\begin{equation*}
\begin{minipage}{0.25\linewidth}
  \begin{center}
    \begin{tikzpicture}[scale=0.72]
      \drawDg{2,3} \setcounter{boxnum}{1}; \foreach \x\y in
      {2/3,1/2,2/2,1/1,2/1} { \node at (\x-0.5,\y-0.5) {\footnotesize
          $\bx{\beta}{\theboxnum}$}; \addtocounter{boxnum}{1}; };
      \node at (1,-.58) {$\beta ^*$};
    \end{tikzpicture}
  \end{center}
  \begin{center}
        \begin{tikzpicture}[scale=0.6]
      \drawDg{2,3} \foreach \x\y in {2/3,2/2} { \node at
        (\x-0.5,\y-0.5) {$1$}; } \node at (0.5,1.5) {$0$};
      \node at (1,-0.62) {\footnotesize \phantom{l}arms\phantom{g}};
    \end{tikzpicture}
    \quad
    \begin{tikzpicture}[scale=0.6]
      \drawDg{2,3} \foreach \x\y in {2/3,1/2} { \node at
        (\x-0.5,\y-0.5) {$0$}; } \node at (1.5,1.5) {$1$};
      \node at (1,-0.62) {\footnotesize legs};
    \end{tikzpicture}
  \end{center}
\end{minipage}
\begin{minipage}{0.75\linewidth}
  \begin{center}
    \begin{tabular}{llll}
    \begin{tikzpicture}[scale=0.4]
      \draw [dashed,gray] (-0.5,1.5)--(3.5,5.5); \draw [dashed,gray]
      (-0.5,0.5)--(3.5,4.5); \draw [dashed,gray]
      (-0.5,-0.5)--(3.5,3.5); \drawskewdg{0/1,0/1}{0}
      \drawskewdg{0/1,0/1,0/1}{2}
      \setcounter{boxnum}{1};
      \foreach \x\y in {3/5,1/2,3/4,1/1,3/3} {
        \node at (\x-0.5,\y-0.5) {\tiny \theboxnum};
        \addtocounter{boxnum}{1};
      }
      \node at (1.5,-1.5) {\(S = \{v_1,v_2,v_3\}\)};
      \node at (3.5,0.5) {\(q^{\shortminus 2}t^4\)};
    \end{tikzpicture}&
    \begin{tikzpicture}[scale=0.4]
      \draw [dashed,gray] (-0.5,1.5)--(2.5,4.5); \draw [dashed,gray]
      (-0.5,0.5)--(3.5,4.5); \draw [dashed,gray]
      (-0.5,-0.5)--(3.5,3.5); \drawskewdg{0/1,0/1}{0}
      \drawskewdg{0/1,-1/1}{2}
      \setcounter{boxnum}{1};
      \foreach \x\y in {2/4,1/2,3/4,1/1,3/3} {
        \node at (\x-0.5,\y-0.5) {\tiny \theboxnum};
        \addtocounter{boxnum}{1};
      }
      \node at (1.5,-1.5) {\(S = \{v_2,v_3\}\)};
      \node at (3.5,0.5) {\(q^{\shortminus 1}t^3\)};
    \end{tikzpicture}&
    \begin{tikzpicture}[scale=0.4]
      \draw [dashed,gray] (-0.5,1.5)--(3.5,5.5); \draw [dashed,gray]
      (-0.5,0.5)--(3.5,4.5); \draw [dashed,gray]
      (-0.5,-0.5)--(4.5,4.5); \drawskewdg{0/1,0/1}{0}
      \drawskewdg{-1/1,-1/0}{3}
      \setcounter{boxnum}{1};
      \foreach \x\y in {3/5,1/2,3/4,1/1,4/4} {
        \node at (\x-0.5,\y-0.5) {\tiny \theboxnum};
        \addtocounter{boxnum}{1};
      }
      \node at (1.5,-1.5) {\(S = \{v_1,v_2\}\)};
      \node at (3.5,0.5) {\(q^{\shortminus 1}t^2\)};
    \end{tikzpicture}&
    \begin{tikzpicture}[scale=0.4]
      \draw [dashed,gray] (-1.5,0.5)--(3.5,5.5); \draw [dashed,gray]
      (-1.5,-0.5)--(3.5,4.5); \draw [dashed,gray]
      (-0.5,-0.5)--(3.5,3.5); \drawskewdg{-1/1}{0}
      \drawskewdg{0/1,0/1,0/1}{2}
      \setcounter{boxnum}{1};
      \foreach \x\y in {3/5,0/1,3/4,1/1,3/3} {
        \node at (\x-0.5,\y-0.5) {\tiny \theboxnum};
        \addtocounter{boxnum}{1};
      }
      \node at (1.5,-1.5) {\(S = \{v_1,v_3\}\)};
      \node at (3.5,0.5) {\(q^{\shortminus 2}t^3\)};
    \end{tikzpicture}\\
    \begin{tikzpicture}[scale=0.4]
      \draw [dashed,gray] (-0.5,1.5)--(2.5,4.5); \draw [dashed,gray]
      (-0.5,0.5)--(3.5,4.5); \draw [dashed,gray]
      (-0.5,-0.5)--(4.5,4.5); \drawskewdg{0/1,0/1}{0}
      \drawskewdg{-2/1}{3}
      \setcounter{boxnum}{1};
      \foreach \x\y in {2/4,1/2,3/4,1/1,4/4} {
        \node at (\x-0.5,\y-0.5) {\tiny \theboxnum};
        \addtocounter{boxnum}{1};
      }
      \node at (1.5,-1.5) {\(S = \{v_2\}\)};
      \node at (3.5,0.5) {\(t\)};
    \end{tikzpicture}&
    \begin{tikzpicture}[scale=0.4]
      \draw [dashed,gray] (-1.5,0.5)--(2.5,4.5); \draw [dashed,gray]
      (-1.5,-0.5)--(3.5,4.5); \draw [dashed,gray]
      (-0.5,-0.5)--(3.5,3.5); \drawskewdg{-1/1}{0}
      \drawskewdg{0/1,-1/1}{2}
      \setcounter{boxnum}{1};
      \foreach \x\y in {2/4,0/1,3/4,1/1,3/3} {
        \node at (\x-0.5,\y-0.5) {\tiny \theboxnum};
        \addtocounter{boxnum}{1};
      }
      \node at (1.5,-1.5) {\(S = \{v_3\}\)};
      \node at (3.5,0.5) {\(q^{\shortminus 1}t^2\)};
    \end{tikzpicture}&
    \begin{tikzpicture}[scale=0.4]
      \draw [dashed,gray] (-1.5,0.5)--(2.5,4.5); \draw [dashed,gray]
      (-1.5,-0.5)--(2.5,3.5); \draw [dashed,gray] (-0.5,-0.5)--(3.5,3.5);
      \drawskewdg{-1/1}{0} \drawskewdg{-1/1,-1/0}{2}
      \setcounter{boxnum}{1};
      \foreach \x\y in {2/4,0/1,2/3,1/1,3/3} {
        \node at (\x-0.5,\y-0.5) {\tiny \theboxnum};
        \addtocounter{boxnum}{1};
      }
      \node at (1.5,-1.5) {\(S = \{v_1\}\)};
      \node at (3.5,0.5) {\(q^{\shortminus 1}t\)};
    \end{tikzpicture}&
    \begin{tikzpicture}[scale=0.4]
      \draw [dashed,gray] (-1.5,0.5)--(1.5,3.5); \draw [dashed,gray]
      (-1.5,-0.5)--(2.5,3.5); \draw [dashed,gray] (-0.5,-0.5)--(3.5,3.5);
      \drawskewdg{-1/1}{0} \drawskewdg{-2/1}{2}
      \setcounter{boxnum}{1};
      \foreach \x\y in {1/3,0/1,2/3,1/1,3/3} {
        \node at (\x-0.5,\y-0.5) {\tiny \theboxnum};
        \addtocounter{boxnum}{1};
      }
      \node at (1.5,-1.5) {\(S = \varnothing\)};
      \node at (3.5,0.5) {\(1\)};
    \end{tikzpicture}
    \end{tabular}
  \end{center}
\end{minipage}
\end{equation*}
\end{example}

\subsection{A formula for \texorpdfstring{$\nabla$}{nabla} on any LLT
polynomial}
\label{ss nabla on LLT}

The operator $\nabla$, introduced in \cite{BeGaHaTe99}, is the linear
operator on symmetric functions which acts diagonally on the basis of
modified Macdonald polynomials $\Htild_{\mu }(X;q,t)$ by $\nabla
\Htild _{\mu } = q^{\mathsf{n}(\mu ^{*})} t^{\mathsf{n}(\mu )}\Htild
_{\mu }$.

In \cite{BHMPS-llt}, we give a raising operator formula for $\nabla$
on any LLT polynomial.  The formula takes a simpler form in
the case that the LLT polynomial is indexed by a tuple of ribbons.  We
state the result for the tuple $\nubold(S)$
in Definition~\ref{d:nuofS}, making use of the
notation $V_\beta, R_{\beta^*}, \widehat{R}_{\beta^*}$ defined in \S\ref{ss
HHL formula}, \eqref{ed:Rbeta}, and \eqref{ed:Rbeta2}.
Also let $A_\beta$ denote the number of attacking pairs in
$\nubold(S)$, which depends only on $\beta$ and not on $S \subseteq
V_\beta$.

\begin{thm}[\cite{BHMPS-llt}]\label{t nabla on LLT}
For $\beta \in \ZZ_{+}^k$ and $S \subseteq V_\beta$, consider the  tuple of ribbons
$\nubold(S)$.
We have the following formula for the
operator $\nabla$ applied to the LLT polynomial
$\Gcal_{\nubold(S)}(X;q)$:

\begin{multline}
\label{eq:nabla on LLT}
\nabla \Gcal_{\nubold(S)}(X;q) = \\
\omega\pol_X \sigmabold
 \bigg( (-qt)^{|V_\beta\setminus S|} q^{A_\beta} \frac{z_1 \cdots z_l \prod_{(\bx{\beta}{i}, \bx{\beta}{j}) \in V_\beta \setminus S} z_i/z_j
\prod_{\alpha_{ij} \in  \widehat{R}_{\beta^*}}
\! \big(1-q\spa  t\, z_i/z_j \big)} { \prod_{\alpha_{ij} \in  R_+} \! \big(1-q \, z_i/z_j\big)
\prod_{\alpha_{ij} \in  R_{\beta^*}} \! \big(1-t \, z_i/z_j\big)}\bigg),
\end{multline}
where $l= |\beta|$ and $R_+ = R_+(\GL _l)$.
\end{thm}

\begin{proof}
Our raising operator formula \cite[Corollary 8.4.1]{BHMPS-llt} for
$\nabla$ on any LLT polynomial reduces to
\eqref{eq:nabla on LLT} using the following facts which hold
when the LLT polynomial is indexed by a tuple of ribbons.

(1) The magic number $p(\nubold(S))$ of $\nubold(S)$ defined in
\cite[\S7.2]{BHMPS-llt} is equal to the number of boxes of
$\nubold(S)$ which are not the first box in a row, which is the same
as $|V_\beta \setminus S|$.

(2) The weight $\lambda$ in \cite[Corollary 8.4.1]{BHMPS-llt}, defined
in \cite[Definition 7.1.2]{BHMPS-llt}, is obtained as follows in the
case that the tuple of skew shapes consists of ribbons: fill the boxes
of each row of $\nubold(S)$ with $1,0,\dots, 0, -1$ or just $0$ for a
row of length 1, and then read this filling by increasing reading
order.  It is then easily seen that $\mathbf{z}^\lambda =
\prod_{(\bx{\beta}{i},\bx{\beta}{j})\in V_\beta \setminus S} z_i/z_j$.

(3) Under the bijection $f \colon \beta \to \nubold(S)$ which takes
the $i$-th box $\bx{\beta}{i}$ of $\dgrm{\beta}$ in reading order to
the $i$-th box of $\nubold(S)$ in reading order, the set
$\{(f(\bx{\beta}{i}),f(\bx{\beta}{j})) : \alpha_{ij} \in
R_{\beta^*}\}$ is exactly the set of non-attacking pairs in
$\nubold(S)$.  Thus $R_{\beta^*}$ agrees with the set of roots denoted
$R_t$ in \cite[Definition 7.1.2 and Remark 7.1.3 (i)]{BHMPS-llt}.
\end{proof}

\subsection{Proof of Corollary \ref{cor:mac formula b 2}}

Applying $\nabla$ to \eqref{eq:HHL formula} and
substituting \eqref{eq:nabla on LLT} into this yields
\begin{align}\label{e proof mac formula}
\nabla \Htild_\mu(X;q,t) = \omega \pol_X \sigmabold \bigg(\frac{z_1
\cdots z_l \cdot \Upsilon \cdot \prod_{\alpha_{ij} \in
\widehat{R}_{\beta^*} } \big(1-q\spa t\, z_i/z_j \big)} {
\prod_{\alpha_{ij} \in R_+} \big(1-q \, z_i/z_j\big)
\prod_{\alpha_{ij} \in R_{\beta^*}} \big(1-t \, z_i/z_j\big)}\bigg),
\end{align}
where
\begin{align}
\Upsilon = \sum_{S \subseteq V_\beta} (-qt)^{|V_\beta \setminus S|}
q^{A_\beta} \prod_{(\bx{\beta}{i}, \bx{\beta}{j})\in S}
q^{-\arm(\bx{\beta}{i})} \spa t^{\leg(\bx{\beta}{i})+1}
\prod_{(\bx{\beta}{i}, \bx{\beta}{j})\in V_\beta \setminus S} z_i/z_j.
\end{align}
Defining $d = \! A_\beta - \mathsf{n}(\mu^*) -
\sum_{(\bx{\beta}{i},\bx{\beta}{j}) \in V_\beta}  \arm(\bx{\beta}{i})$
and using $\mathsf{n}(\mu) = \! \sum_{(\bx{\beta}{i},
\bx{\beta}{j})\in V_\beta}  (\leg(\bx{\beta}{i})+1)$,
the quantity $\Upsilon$ can be
simplified as follows:
\begin{align*}
\Upsilon &= q^{\mathsf{n}(\mu^*)+d} \, t^{\mathsf{n}(\mu)} \sum_{S
\subseteq V_\beta}\Big( (-qt)^{|V_\beta\setminus S|}
\prod_{(\bx{\beta}{i}, \bx{\beta}{j})\in V_\beta}
\!\!\!\!  q^{\arm(\bx{\beta}{i})} \, t^{-\leg(\bx{\beta}{i})-1} \\
 & \hphantom{q^{\mathsf{n}(\mu^*)+d} \, t^{\mathsf{n}(\mu)} \sum_{S
\subseteq V_\beta}\Big( (-qt)^{|V_\beta\setminus S|}
\prod_{(\bx{\beta}{i}, \bx{\beta}{j})\in V_\beta}}
 \times \prod_{(\bx{\beta}{i}, \bx{\beta}{j})\in S} \!\!\!\!
q^{-\arm(\bx{\beta}{i})} \, t^{\leg(\bx{\beta}{i})+1} \!\!\!
\prod_{(\bx{\beta}{i}, \bx{\beta}{j})\in V_\beta \setminus S}
\!\!\!\!\!\! z_i/z_j \Big)
\\
 &= q^{\mathsf{n}(\mu^*)+d} \spa t^{\mathsf{n}(\mu)}
\sum_{S \subseteq V_\beta} \prod_{(\bx{\beta}{i}, \bx{\beta}{j})\in
V_\beta\setminus S} \!\! \big(-q^{\arm(\bx{\beta}{i})+1} \,
t^{-\leg(\bx{\beta}{i})} z_i/z_j \big)\\
 &= q^{\mathsf{n}(\mu^*)+d} \spa t^{\mathsf{n}(\mu)}
\prod_{(\bx{\beta}{i}, \bx{\beta}{j}) \in V_\beta} \big(1-
q^{\arm(\bx{\beta}{i})+1}\spa t^{-\leg(\bx{\beta}{i})} z_i/z_j \big).
\end{align*}
Thus, plugging this back in for $\Upsilon$ in \eqref{e proof mac
formula} and recalling the definition of $\HS_{\beta^*}(\zz;q, t)$
from \eqref{et:mac formula b}, we have
\[ \nabla \Htild_\mu(X;q,t) =
q^{\mathsf{n}(\mu^*)+d} \spa t^{\mathsf{n}(\mu)}
\omega\pol_X \! \big(z_1\cdots z_l \, \HS_{\beta^*}(\zz;q, t) \big).
\]
By the definition of $\nabla$, this implies $\Htild_\mu(X;q,t) = q^{d}
\omega \pol_X \! \big(z_1\cdots z_l \, \HS_{\beta^*}(\zz;q, t) \big).$
It remains to show that $d=0$.  This follows from the fact that the
coefficient of $s_{1^l}(X) = \pol_X(z_1\cdots z_l)$ in the Schur
expansion of both $\omega\Htild_\mu(X;q,t)$ and $\pol_X \!
\big(z_1\cdots z_l \, \HS_{\beta^*}(\zz;q, t) \big)$ is 1; the former
is well known, while the latter can be seen directly by expanding the
series $z_1\cdots z_l \, \HS_{\beta^*}(\zz;q, t)$ to see that it is
equal to $\sigmabold(z_1\cdots z_l)= \chi_{1^l}$ plus terms of the
form $a_\nu \chi_\nu$ for $\nu > 1^l$ in dominance order.

\section{The \texorpdfstring{$m,n$}{m,n}-Macdonald polynomials}
\label{s:m,n-Macdonalds}
\setcounter{subsection}{1}

For every pair of coprime positive integers $(m,n)$, the action
of the Burban-Schiffmann elliptic Hall algebra $\Ecal $ on $\Lambda
(X)$ gives rise to a family of symmetric functions that we call
$m,n$-Macdonald polynomials.  The subfamily of $1,n$-Macdonald polynomials is
closely connected with the Macdonald series $\Hbold _{\mu }$
from Definition \ref{d mac formula}.
In this section we
will construct raising operator formulas for all $m,n$-Macdonald
polynomials, reducing to
Corollary \ref{cor:mac formula b 2}
in the case
$m=n=1$.

To define $m,n$-Macdonald polynomials we need to recall some
facts about the algebra $\Ecal $, defined by Burban and Schiffmann
\cite{BurbSchi12} in terms of Hall algebras of coherent sheaves on
elliptic curves.  For each pair of coprime integers $(m,n)$, the
algebra $\Ecal $ contains a distinguished subalgebra $\Lambda
(X^{m,n}) $ isomorphic to the algebra of symmetric functions over $\kk
$; these subalgebras generate $\Ecal $, subject to relations given in
\cite{BurbSchi12}.  There is also an action of $\Ecal $ on the space
of symmetric functions $\Lambda (X)$, constructed by Schiffmann and
Vasserot \cite{SchiVass13}.  Our notation here is the same as in
\cite{BHMPS-lw,BHMPS-llt,BHMPS-delta,BHMPS-paths}---in particular, we
we use the version of the action of $\Ecal$ on $\Lambda(X)$ given by
\cite[Proposition~3.3.1]{BHMPS-paths}.  The translation between our
notation and that of \cite{BurbSchi12,SchiVass13} can be found in
\cite[\S \S3.2--3.3]{BHMPS-paths}; the defining relations of $\Ecal $
written in our notation are given in \cite[\S 3.2]{BHMPS-delta}.

\begin{defn}
\label{d mn mac} Set $M = (1-q)(1-t) \in \kk$.  For coprime positive
integers $m,n$, define the \emph{$m,n$-Macdonald polynomial}
$\Htild_{\mu}^{m,n} = \Htild_{\mu}^{m,n}(X;q,t)$ by
\begin{align}\label{e:mn-mac}
\Htild_{\mu}^{m,n} =  \Htild_\mu[-M X^{m,n}] \cdot 1,
\end{align}
the symmetric function obtained by
acting on $1 \in \Lambda(X)$ with
(a plethystic transformation of) a modified Macdonald polynomial in
the distinguished subalgebra $\Lambda(X^{m,n}) \subseteq \Ecal$.
\end{defn}

\begin{remark}
(i) For context, we note that
$e_k[-M X^{m,n}] \cdot 1$
defines the symmetric function side
of the $(km,kn)$-shuffle theorem of \cite{BeGaSeXi16,Mellit16}.

(ii) $\Htild_{\mu}^{m,n}$ is a homogeneous symmetric function of degree $n|\mu|$, as follows from the definition of the action of the Schiffmann algebra on symmetric functions \cite[Proposition 3.3.1]{BHMPS-paths}.
\end{remark}

\subsection{The \texorpdfstring{$1,n$}{1,n}-Macdonald polynomials}
\label{ss:m=1-and-n=1}

By \cite[Lemma 3.5.1]{BHMPS-llt}, $\Htild_\mu^{m,1} = \Htild_\mu[-M
X^{m,1}] \cdot 1 = \nabla^m \Htild_\mu = q^{m \spa \mathsf{n}(\mu^*)}t^{m \spa
\mathsf{n}(\mu)} \Htild_\mu$, so this case is familiar.  The $1,n$-Macdonald
polynomials, on the other hand, carry essentially the same data as the
series $\HS_\mu(\zz;q,t)$ by the following theorem,
which will be proven as part of a more general result below (see
Remark \ref{r mac formula mn} (ii)).

\begin{thm}\label{t mac formula 1n}
For any partition $\mu$ of $l$ and any rearrangement $\beta$ of
$\mu^{*}$,
\begin{align}\label{et:mac formula 1n}
\Htild_{\mu}^{1,n}(X;q,t) = \omega\pol_X \!
\Big(q^{\mathsf{n}(\mu^*)}t^{\mathsf{n}(\mu)}\, (z_1\cdots z_l)^n \,
\HS_{\beta^*}(\zz;q, t) \Big).
\end{align}
Hence also $\HS_{\beta^*}(\zz;q, t) = q^{-\mathsf{n}(\mu^*)}t^{-\mathsf{n}(\mu)} \lim_{n \to
\infty} (z_1\cdots z_l)^{-n}(\omega\Htild_{\mu}^{1,n})(z_1,\dots, z_l).$
\end{thm}

This also shows that $\HS_{\beta^*}(\zz;q, t)$ depends only on the partition rearrangement $\mu^*$ of $\beta$,
thus establishing Theorem \ref{t mac formula 2}.
Theorem \ref{t mac formula 1n} also shows that
the following is equivalent to Conjecture~\ref{conj:Hmu-positivity}.

\begin{conj}
\label{cj 1n Schur positive}
The $1,n$-Macdonald polynomials $\Htild_{\mu}^{1,n}(X;q,t)$ are Schur
positive.
\end{conj}
\begin{remark}
For $m \ne 1$ and $n \ne 1$, the $m,n$-Macdonald polynomials are typically not Schur positive.
\end{remark}

\begin{example}
\label{ex 1nmacs}
The Schur expansions of the  $1,n$-Macdonald polynomials   $\Htild^{1,n}_\mu(X;q,t)$ for  $n=2$ and $|\mu| = 2, 3$, written as in \eqref{et:mac formula 1n}, are
\begin{align*}
\Htild^{1,2}_{2} \, &=  \omega \, q \big( q^2s_{4} + q \spa s_{31} + s_{22} \big)\\
\Htild^{1,2}_{11} \, &= \omega \,t \big( t^2s_{4} + t \spa s_{31} + s_{22} \big) \\
\Htild^{1,2}_{3} \, &=  \omega \,q^3 \big( q^6s_{6} + (q^5 + q^4)s_{51} + (q^4 + q^3 + q^2)s_{42} +
q^3s_{33} + q^3s_{411} + (q^2 + q)s_{321} + s_{222}\big) \\
\Htild^{1,2}_{21} \, &=  \omega \, qt \big( q^2t^2s_{6} + (q^2t + qt^2)s_{51} + (q^2 + qt + t^2)s_{42} +
qt \spa s_{33} + qt \spa s_{411} + (q + t)s_{321} + s_{222} \big) \\
\Htild^{1,2}_{111} \, &= \omega\, t^3 \big( t^6s_{6} + (t^5 + t^4)s_{51} + (t^4 + t^3 + t^2)s_{42} +
t^3s_{33} + t^3s_{411} + (t^2 + t)s_{321} + s_{222} \big)
\end{align*}
\end{example}

\begin{prop}
\label{prop 1,n mac q,t symmetry}
The $1,n$-Macdonald polynomials satisfy the same $q,t$ symmetry property as ordinary modified Macdonald polynomials:
$\Htild^{1,n}_\mu(X;q,t) = \Htild^{1,n}_{\mu^*}(X;t,q)$.
\end{prop}
\begin{proof}
This follows from the symmetry property for Macdonald polynomials, $\Htild_\mu(X;q,t) = \Htild_{\mu^*}(X;t,q)$, and the fact that for any symmetric function $f$ with coefficients in $\QQ$, $f[-M X^{m,n}] \cdot 1$ is symmetric in $q$ and $t$, which in turn follows from the description of the action of the Schiffmann algebra on symmetric functions in \cite[Proposition 3.3.1]{BHMPS-paths}.
\end{proof}

We obtain several specializations of the $1,n$-Macdonald polynomials
easily from the raising operator formula in Theorem \ref{t mac formula
1n}.

It is convenient to work with another variant of the Hall-Littlewood
polynomials $H_\mu(X;t) \defeq
t^{\mathsf{n}(\mu)}\Htild_\mu(X;0,t^{-1})$.

\begin{prop}
\label{p specialize 1n} Let $\mu$ be a partition of $l$ with transpose
$\nu = \mu^*$.  The $q=t=1$, $q=1$, and $q=0$ specializations of the
$1,n$-Macdonald polynomials are given by
\begin{gather}
\label{e 1nmac q1t1}
\Htild^{1,n}_\mu(X;1,1) = e_{(n^l)}(X), \\
\label{e 1nmac q1}
\Htild^{1,n}_\mu(X;1,t) = t^{\mathsf{n}(\mu)} \prod_{r = 1}^{\mu_1}\omega H_{(n^{\nu_r})}(X;t), \\
\label{e 1nmac q0}
\big(q^{-\mathsf{n}(\mu^*)} t^{-\mathsf{n}(\mu)}\Htild^{1,n}_\mu(X;q,t) \big)|_{q=0}
= \omega \pol_X \! \Big( \sigmabold\Big(\frac{(z_1\cdots z_l)^n}{\prod_{\alpha_{ij} \in B_\mu}
(1 - t \spa z_i/z_j)}\Big) \Big),
\end{gather}
where $B_\mu$ is as defined before \eqref{e HL formula from 1s}.
\end{prop}

\begin{remark}
Like the familiar right hand sides of \eqref{e 1nmac q1t1} and
\eqref{e 1nmac q1}, the right hand side of \eqref{e 1nmac q0} is also
well studied.  Its Schur expansion coefficients are instances of the
generalized Kostka polynomials introduced by Shimozono-Weyman
\cite{ShimWeym00}
corresponding to the
sequence of rectangle shapes
 $(n^{\mu_\ell}), \dots, (n^{\mu_2}), (n^{\mu_1})$ for $\ell =
\ell(\mu)$.
\end{remark}

\begin{proof}
By the same argument as in \S\ref{ss specialization}, the $q=t=1$
specialization of $\pol_X ((z_1\cdots z_l)^n\HS_\mu)$ is the complete
homogeneous symmetric function $h_{(n^l)}$.  Hence \eqref{e 1nmac
q1t1} follows from Theorem \ref{t mac formula 1n}.  Similarly,
\eqref{e 1nmac q0} follows from Theorem \ref{t mac formula 1n} and
from noting that $\pol_X ((z_1\cdots z_l)^n\HS_\mu(\zz; 0,t))$ can be
simplified just as $\pol_X (z_1\cdots z_l \, \HS_\mu(\zz; 0,t))$ was
in \S\ref{ss specialization}.

We now prove \eqref{e 1nmac q1}. By Theorem \ref{t mac formula 1n} and
\eqref{e:Weyl-symmetrizer 2},
\begin{align}\label{e 1n mac start}
q^{-\mathsf{n}(\mu^*)}t^{-\mathsf{n}(\mu)} \omega\spa
\Htild^{1,n}_{\mu }(X;q,t) = \mathbf{h}_X\Big( \spa (z_1\cdots z_l)^n
\prod_{\alpha_{ij} \in R_+}(1-z_i/z_j)\widehat{\HS}_{\mu}\Big),
\end{align}
where
$\widehat{\HS}_\mu(\zz;q,t)$ is the expression inside the
$\sigmabold(\cdot)$ on the right side of \eqref{et:mac formula a}.
Let $C_1, \dots, C_{\mu_1}$ denote the columns of the diagram of $\mu$
and note that $R_\mu \setminus \widehat{R}_\mu$ is equal to
$\bigsqcup_{r = 1}^{\mu_1} \{\alpha_{ij} :
\bx{\mu}{j}=\south(\bx{\mu}{i}), \, \bx{\mu}{i} \in C_r \}$.  Hence
setting $q=1$ in \eqref{e 1n mac start} yields
\begin{align}
t^{-\mathsf{n}(\mu)} \omega \Htild^{1,n}_{\mu }(X;1,t)
& = \mathbf{h}_X\bigg(
\prod_{r = 1}^{\mu_1} \frac{\big(\prod_{\bx{\mu}{i} \in C_r} z_i^n\big) \prod_{\alpha_{ij} \in R_\mu \setminus \widehat{R}_\mu, \, \bx{\mu}{i} \in C_r }
\big(1-  t^{-\leg(\bx{\mu}{i})} z_i/z_j \big)}
{\prod_{\alpha_{ij} \in R_\mu \setminus \widehat{R}_\mu, \, \bx{\mu}{i} \in C_r}\big(1-t z_i/z_j \big)}
\bigg) \\
\label{e 1n mac step2}
&  =
\prod_{r = 1}^{\mu_1}\mathbf{h}_X\bigg(
\frac{\big(\prod_{\bx{\mu}{i} \in C_r} z_i^n\big) \prod_{\alpha_{ij} \in R_\mu \setminus \widehat{R}_\mu, \, \bx{\mu}{i} \in C_r }
\big(1-  t^{-\leg(\bx{\mu}{i})} z_i/z_j \big)}
{\prod_{\alpha_{ij} \in R_\mu \setminus \widehat{R}_\mu\, \bx{\mu}{i} \in C_r}\big(1-t z_i/z_j \big)}
\bigg),
\end{align}
where the second equality follows from \eqref{e hX equals h} and the
fact that $h_\gamma = h_\delta$ for $\gamma, \delta \in \ZZ^l$ which
are rearrangements of each other.

The factor $\mathbf{h}_X(\cdot)$ in \eqref{e 1n mac step2} for a given
index $r$ is equal to $t^{-\mathsf{n}(1^{\nu_r})} \omega
\Htild^{1,n}_{(1^{\nu_r})}(X;1,t)$, by the computation we have just
done, but with the partition $1^{\nu_r}$ in place of $\mu $.  It
follows that
\begin{align}
\label{e 1n mac factor}
\Htild^{1,n}_{\mu }(X;1,t) =
\prod_{r = 1}^{\mu_1}\Htild^{1,n}_{(1^{\nu_r})}(X;1,t),
\end{align}
Finally, we will show that each $\Htild^{1,n}_{(1^{\nu_r})}(X;1,t)$ is
essentially a Hall-Littlewood polynomial.
Using the particularly simple form of the series $\HS_{\mu}$ when $\mu
= (d)$ is a single row,
Theorem \ref{t mac formula 1n} and \eqref{e HL formula from mu} yield
$q^{-\binom{d}{2}}\omega \spa \Htild^{1,n}_{(d) }(X;q,t) =
H_{(n^d)}(X;q)$.
By Proposition~\ref{prop 1,n mac q,t symmetry}, $\Htild^{1,n}_{\mu }(X;q,t) =
\Htild^{1,n}_{\mu^* }(X;t,q)$, and thus
\begin{align}
\label{e 1n mac 1d}
t^{-\binom{d}{2}}\omega \spa \Htild^{1,n}_{(1^d) }(X;q,t) = H_{(n^d)}(X;t).
\end{align}
Formula \eqref{e 1nmac q1} now follows from \eqref{e 1n mac factor}
and \eqref{e 1n mac 1d}.
\end{proof}

\subsection{A raising operator formula for the
\texorpdfstring{$m,n$}{m,n}-Macdonald polynomials}
\label{ss:m,n-Mac-formula}

Our raising operator formula for $\Htild_\mu$ is a special case of
a raising operator formula for $\Htild_\mu^{m,n}$ which can be
obtained in a similar way.  Specifically, we
combine the Haglund-Haiman-Loehr formula
(Theorem \ref{t HHL formula}) with an $m,n$ version of Theorem \ref{t
nabla on LLT} from \cite{BHMPS-llt}.  This latter result requires some
notation involving the dilation of a column diagram $\beta$.

For  $\beta \in \ZZ_{+}^k$, let  $m \beta = (m \beta_1,\dots, m \beta_k)$. We think of the column diagram  $\dgrm{(m\beta)}$ of  $m \beta$
as the result of dilating the column diagram $\beta^*$ of  $\beta$ vertically by a factor of $m$,
so that each box of  $\dgrm{\beta}$ gives rise to  $m$ boxes of  $\dgrm{(m \beta)}$. To be more precise, define the map of boxes
\begin{align}
\tau \colon \dgrm{(m \beta)} \to \dgrm{\beta},  \ \ \ \  (i,j) \mapsto  (i,\lfloor (j-1)/m \rfloor + 1).
\end{align}
Thus, in the dilation process, each box  $b$ of  $\dgrm{\beta}$ gives rise to a set $\tau^{-1}(b)$ of $m$ boxes of $\dgrm{(m \beta)}$, called a \emph{dilated box}, 
which forms a contiguous subset of a column of  $\dgrm{(m \beta)}$.

As in \S\ref{ss:column-diagram-formula}, we write
$\bx{\beta}{1},\dots, \bx{\beta}{d}$ for the boxes of $\dgrm{\beta}$
listed increasing reading order and $\bx{m\beta}{1},\dots,
\bx{m\beta}{l}$ for the boxes of $\dgrm{(m\beta)}$ in increasing
reading order, where $d = |\beta|$ and $l = dm$.  Define a map $\str
\colon V_\beta \to V_{m \beta}$ which takes a vertical domino
$(\bx{\beta}{i}, \bx{\beta}{j})$ of $\dgrm{\beta}$ to the vertical
domino $(\bx{m\beta}{i'}, \bx{m\beta}{j'})$ of $\dgrm{(m \beta)}$ such
that $\tau(\bx{m\beta}{i'}) = \bx{\beta}{i}$ and
$\tau(\bx{m\beta}{j'}) = \bx{\beta}{j}$.  Note that $\str(V_\beta)$ is
the set of vertical dominoes of the dilated diagram $\dgrm{(m \beta)}$
which straddle two dilated boxes.

\begin{example}
For \(m=3\) and \(\beta = (2,1,2)\), the dilated boxes of
$\dgrm{(m\beta)}$ are shown below along with the labeling of the boxes
of $\dgrm{\beta}$ and $\dgrm{(m\beta)}$ by reading order.  To clarify
the definitions of $\tau$ and $\str$, note that
\(\tau^{-1}(\bx{\beta}{3}) = \{\bx{m\beta}{7}, \bx{m\beta}{10},
\bx{m\beta}{13}\}\) and \(\str(V_\beta) =
\{(\bx{m\beta}{5},\bx{m\beta}{7}),
(\bx{m\beta}{6},\bx{m\beta}{9})\}\).
  \begin{equation*}
      \begin{tikzpicture}[xscale = 1.05, yscale=0.8]
    \drawDg{2,1,2}
    \setcounter{boxnum}{1};
    \foreach \x\y in {1/2,3/2,1/1,2/1,3/1} {
      \node at (\x-0.5,\y-0.5) {\tiny $\bx{\beta}{\theboxnum}$};
      \addtocounter{boxnum}{1};
     };
  \node at (1.5,-0.7) {$\dgrm{(\beta)}$};
  \end{tikzpicture}
  \raisebox{2cm}{\ \ \ {\Large \text{\(\overset{\, \tau}{\leftarrow}\)}} \  }
  \begin{tikzpicture}[xscale = 1.05, yscale=0.8]
    \drawDg{6,3,6}
    \setcounter{boxnum}{1};
    \foreach \x\y in {1/6,3/6,1/5,3/5,1/4,3/4,1/3,2/3,3/3,1/2,2/2,3/2,1/1,2/1,3/1} {
      \node at (\x-0.5,\y-0.5) {\tiny $\bx{m\beta}{\theboxnum}$};
      \addtocounter{boxnum}{1};
     };
  \node at (1.5,-0.7) {$\dgrm{(m\beta)}$};
  \end{tikzpicture}
  \hspace{2cm}
  \begin{tikzpicture}[scale=0.8]
    \foreach \x\y in {1/1,2/1,3/1,1/4,3/4} {
      \draw[draw = none, fill=blue!29] (\x-1+.27,\y-1+.31) rectangle
      (\x-.27, \y+2-.31);
    }
    \drawDg{6,3,6}
  \node at (1.5,-0.7) {Dilated boxes of $\dgrm{(m\beta)}$};
  \end{tikzpicture}
  \end{equation*}
\end{example}

Given $(m,n)\in \ZZ _{+}\times \ZZ_+ $, we define the sequence of $m$
integers as in \cite[(104)]{BHMPS-llt}
\begin{equation}\label{e:bbold}
\bb (m,n)_{i} = \lceil i n/m \rceil -\lceil (i-1) n/m \rceil \qquad
(i=1,\ldots,m).
\end{equation}
We then define, for any $\beta \in \ZZ_{+}^k$ of size $d= |\beta|$, a
weight $\bb(m,n,\beta) \in \ZZ^{dm}$ as follows: fill each dilated box
of $\dgrm{(m \beta)}$ with the sequence $\bb(m,n)$ from north to
south, and then read this filling by the reading order on $\dgrm{(m \beta)}$.
Equivalently, $\bb(m,n,\beta)_r = \bb(m,n)_a$, where $a$ is the
integer $a\in \{1,\dots,m\}$ such that $a\equiv -j+1 \pmod{m}$ for $j$
the row index of the $r$-th box of $\dgrm{(m\beta)}$ in reading order (i.e.
$\bx{m \beta}{\, r} = (i,j)$ for some $i$).

\begin{thm}[\cite{BHMPS-llt}]\label{t LLT mn}
Let $m,n$ be coprime positive integers, let $\beta \in
\ZZ_{+}^k$, and set $d = |\beta|$, $l = dm$.
For $S \subseteq
V_\beta$, let $\nubold(S)$ be the  $k$-tuple of ribbons in
Definition~\ref{d:nuofS}.  The action of the LLT polynomial
$\Gcal_{\nubold(S)}[-MX^{m,n}] \in \Lambda(X^{m,n})$ on 1 is given by
\begin{multline}
\label{eq:LLT mn}
\Gcal_{\nubold(S)}[-MX^{m,n}] \cdot 1 =
 \omega \pol_X \sigmabold
\bigg( (-qt)^{|V_\beta\setminus S|} q^{A_\beta+(m-1)\mathsf{n}(\beta^+)}
 \mathbf{z}^{\bb(m,n,\beta)} \\
\times \frac{ \prod_{(\bx{m\beta }{i}, \bx{m\beta }{j})
\in \str(V_\beta \setminus S)} z_i/z_j \prod_{\alpha_{ij} \in
\widehat{R}_{(m\beta) ^*}} \! \big(1-q\spa t\, z_i/z_j \big)} {
\prod_{\alpha_{ij} \in R_{+}(\GL_{l})} \! \big(1-q \, z_i/z_j\big)
\prod_{\alpha_{ij} \in R_{(m\beta )^*}} \! \big(1-t \,
z_i/z_j\big)}\bigg),
\end{multline}
where
$A_\beta$ is the number of attacking pairs of  $\nubold(S)$ as in \S\ref{ss nabla on LLT},
$\beta^+$ is the partition rearrangement of $\beta$, and
$R_{(m\beta)^*}, \spa \widehat{R}_{(m\beta)^*} \subseteq
R_+(\GL_l)$ are as in (\ref{ed:Rbeta}, \ref{ed:Rbeta2}).
\end{thm}

\begin{proof}
This is obtained by combining \cite[Theorem~8.3.1]{BHMPS-llt} and
\cite[Proposition 2.3.2]{BHMPS-lw} and specializing to the case that
the LLT polynomial is indexed by the tuple of ribbons $\nubold(S)$.
The notation here and that in \cite[Theorem~8.3.1]{BHMPS-llt} are
matched using the discussion in the proof of Theorem \ref{t nabla on
LLT} and the following: the weight $\bb(m,n,\beta)$ is the same as
$\tilde{\bb}$ defined in \cite[(106)]{BHMPS-llt}, and the weight
$\lambda$ defined in \cite[(105)]{BHMPS-llt} satisfies $\zz^\lambda =
\mathbf{z}^{\bb(m,n,\beta)} \prod_{(\bx{m\beta}{i},
\bx{m\beta}{j}) \in \str(V_\beta \setminus S)} z_i/z_j$.
\end{proof}

\begin{thm}\label{t mac formula mn}
Let $m,n$ be coprime positive integers.  Given $\beta \in \ZZ_{+}^k$
and setting $d = |\beta|$ and $l = dm$,  define the \emph{$m,n$-Macdonald
series} by
\begin{multline}
\label{et:mac formula mn a}
\HS_{\beta^*}^{m,n}(\zz;q,t) = \\
\ \sigmabold
\Bigg(
\mathbf{z}^{\bb(m,n,\beta)} \
\frac{
\prod\limits_{
(\bx{m \beta}{i},\bx{m \beta}{j}) \in \str(V_\beta)
}
\!\big(1- q^{\arm(\tau(\bx{m\beta}{i}))+1}\spa t^{-\leg(\tau(\bx{m\beta}{i}))} z_i/z_j \big)
\prod\limits_{\alpha_{ij} \in \widehat{R}_{(m\beta)^*}}\!\!
\big(1-q\spa  t\, z_i/z_j \big)}
 {\prod_{\alpha_{ij} \in R_+(\GL_l)} \!\big(1-q \, z_i/z_j\big)
\prod_{\alpha_{ij} \in R_{(m\beta)^*}}\!\big(1-t \, z_i/z_j\big)}\Bigg),
\end{multline}
which we regard as an infinite formal linear combination of
irreducible $\GL _{l}$ characters using the convention of Remark
\ref{rem:geometric-series}.  Then, for any partition $\mu$ and any
rearrangement $\beta$ of $\mu^{*}$,
\begin{align}
\label{et:mac formula mn b}
\Htild_{\mu}^{m,n}(X;q,t)
 = \omega \pol_X \! \Big( q^{m \spa \mathsf{n}(\mu^*)} \spa t^{\mathsf{n}(\mu)} \spa \HS_{\beta^*}^{m,n}(\zz;q, t) \Big).
\end{align}
\end{thm}

\begin{remark}
\label{r mac formula mn}
(i) In \eqref{et:mac formula mn a}, the indices of $z_i$ correspond to the
boxes of the dilated diagram $\dgrm{(m \beta)}$, while the arm and leg are
taken with respect to the original diagram $\dgrm{\beta}$.

(ii) Since $\bb(1,n,\beta) = n^l$, $\HS^{1,n}_{\beta ^{*}} = (z_1\cdots
z_l)^n\HS_{\beta^*}$, where $\HS_{\beta^*}$ is as in \eqref{et:mac
formula a 2}.  Hence Theorem \ref{t mac formula mn} specializes to
Corollary \ref{cor:mac formula b 2}
when $(m,n) = (1,1)$,
and also proves Theorem \ref{t mac formula 1n} by setting $m=1$
with $n$ arbitrary.

(iii) Expanding on (ii), the series $\HS^{m,n}_{\beta ^{*}}$
simultaneously encodes the $m,n$-Macdonald polynomials $\{H^{m,n'}_\mu
: n' \in (n+m\ZZ) \cap \ZZ_+\}$, in the following sense:
\begin{align}
\label{e stable shift}
\Htild_{\mu}^{m,n+am}
 = \omega \pol_X \! \Big( q^{m \spa \mathsf{n}(\mu^*)} \spa t^{\mathsf{n}(\mu)} \spa (z_1\cdots z_l)^a \spa \HS_{\beta^*}^{m,n}\Big).
\end{align}
This follows from Theorem \ref{t mac formula mn}
using that $\bb(m,n+a\spa m,\beta) = \bb(m,n,\beta) + a^l$, which in turn holds by $\bb(m,n+a\spa m) = \bb(m,n)+a^m$.
\end{remark}

Theorem \ref{t mac formula mn} and Remark \ref{r mac formula mn} (iii) have the following corollary.

\begin{cor}
The $m,n$-Macdonald series $\HS^{m,n}_{\beta^*}$ depends only on the
multiset of parts of $\beta$ and not on their order:
$\HS^{m,n}_{\beta^*} = \HS^{m,n}_{\gamma^*}$ for any rearrangement
$\gamma$ of $\beta$.
\end{cor}

\begin{proof}[Proof of Theorem \ref{t mac formula mn}] By a plethystic
transformation and change of variables we can replace $X$ in
\eqref{eq:HHL formula} with $-MX^{m,n}$.  Letting both sides act on $1
\in \Lambda(X)$ and then substituting in \eqref{eq:LLT mn} yields
\begin{align}
\label{e mac formula mn step} & \Htild_\mu[-M X^{m,n}] \cdot 1 =
\omega \pol_X \sigmabold \bigg( \frac{\mathbf{z}^{\bb(m,n,\beta)}
\cdot \Upsilon \cdot \prod_{\alpha_{ij} \in
\widehat{R}_{(m\beta)^*} } \big(1-q\spa t\, z_i/z_j \big)}
{\prod_{\alpha_{ij} \in R_+} \big(1-q \, z_i/z_j\big)
\prod_{\alpha_{ij} \in R_{(m\beta)^*}} \big(1-t \,
z_i/z_j\big)}\bigg),
\end{align}
where
\begin{align}
\Upsilon = \sum_{S \subseteq V_\beta} (-qt)^{|V_\beta \setminus S|}
q^{A_\beta+(m-1)\mathsf{n}(\mu^*)} \prod_{(\bx{\beta}{i},
\bx{\beta}{j})\in S} \!\!\! q^{-\arm(\bx{\beta}{i})} \spa
t^{\leg(\bx{\beta}{i})+1} \prod_{(\bx{m\beta}{i},
\bx{m\beta}{j})\in \str(V_\beta \setminus S)} \!\!\!\!z_i/z_j.
\end{align}
The proof of Corollary \ref{cor:mac formula b 2} establishes that $A_\beta =
\mathsf{n}(\mu^*)+ \sum_{(\bx{\beta}{i},\bx{\beta}{j}) \in V_\beta}
\arm(\bx{\beta}{i})$, as can also be checked combinatorially.  Using this, $\Upsilon$ simplifies essentially
the same way it did in that proof:
\begin{align}
\Upsilon \ =\ q^{m \spa \mathsf{n}(\mu^*)} \spa t^{\mathsf{n}(\mu)}
\prod_{(\bx{m\beta}{i}, \bx{m\beta}{j}) \in
\str(V_\beta)} \big(1- q^{\arm(\tau(\bx{m\beta}{i}))+1}\spa
t^{-\leg(\tau(\bx{m\beta}{i}))} z_i/z_j \big).
\end{align}
Plugging this back into \eqref{e mac formula mn step} completes the proof.
\end{proof}

We obtain yet other expressions for the modified Macdonald polynomials
from Theorem \ref{t mac formula mn}.
\begin{cor}
\label{cor:m1 mac formula} Let  $\mu$ be a partition of  $d$ and
$\beta$ any rearrangement of  $\mu^*$. The modified Macdonald polynomial
$\Htild_{\mu }(X;q,t)$ can be expressed in terms of the $m,1$-Macdonald
series as follows:
\begin{align}
\label{ec:m1 mac formula} \Htild_\mu(X;q,t) = \omega \pol_X \! \big(
t^{(-m+1)\mathsf{n}(\mu)}\,  \HS_{\beta^*}^{m,1}(\zz;q, t)\big).
\end{align}
\end{cor}
\begin{proof}
By \cite[Lemma 3.5.1]{BHMPS-llt},
$\Htild_{\mu}[-MX^{m,1}] \cdot 1 = \nabla^m \Htild_\mu$.  Substitute this into the left side of \eqref{et:mac formula mn b}.
\end{proof}

\section{A raising operator formula for the Macdonald polynomials
\texorpdfstring{$Q_\mu(X;q,t)$}{Q mu(X;q,t)} and
\texorpdfstring{$J_\mu(X;q,t)$}{J mu(X;q,t)}}
\label{s integral form J}
\setcounter{subsection}{1}

Our formulas for the modified Macdonald polynomials
$\Htild_\mu(X;q,t)$ can be converted into formulas for the integral
form Macdonald polynomials $J_{\mu }(X;q,t)$ and for the
$Q_\mu(X;q,t)$ which differ from the $J_\mu$ by a scalar factor.
Recall that $\Htild_{\mu }(X;q,t) =
t^{\mathsf{n}(\mu )} J_{\mu}[\frac{X}{1-t^{-1}};q,t^{-1}]$ as in \eqref{e:H-tilde}, hence $\Htild_{\mu }[X(1-t^{-1});q,t] =
t^{\mathsf{n}(\mu )} J_{\mu}(X;q,t^{-1})$, or equivalently
\begin{align}\label{e J in terms of Htild}
J_{\mu}(X;q,t) = t^{\mathsf{n}(\mu )} \Htild _{\mu }[X(1-t);q,t^{-1}].
\end{align}

Let $\mathbf{e}'_X \colon \kk[z_1^{\pm1},\dots, z_l^{\pm1}] \to
\Lambda$ denote the linear operator determined by
\begin{align}\label{e:e'}
\mathbf{e}'_X(\zz^\gamma) &= e_\gamma[X(1-t)] = \prod_{i = 1}^l
\sum_{j=0}^{\gamma _i} (-t)^{j} e_{\gamma _i -j}(X) h_{j}(X),
\end{align}
where for  $\gamma \in \ZZ^l$, we define $e_\gamma =
e_{\gamma_1}\cdots e_{\gamma_l}$ to be the product of elementary symmetric functions,
with $e_d$ for $d\leq 0$ interpreted as $e_0 =1$, or $e_d = 0$ for $d< 0$.
We extend  $\mathbf{e}'_X$ to an operator on  formal linear combinations of monomials
just as we did for  $\sigmabold$ and  $\mathbf{h}_X$ in \S\ref{ss:Weyl-etc}.
The operators $\mathbf{e}'_X$ and  $\mathbf{h}_X$ are related as follows: for a formal linear combination
$g = \sum_{\gamma \in \ZZ^l} c_\gamma \zz^\gamma$
such that  $\mathbf{h}_X(g)$ is a symmetric function and not just
an infinite formal sum of symmetric functions,
 $\mathbf{e}'_X(g) = \theta \circ \omega \circ \mathbf{h}_X(g)$, where
$\theta \colon \Lambda \to \Lambda$ is the automorphism sending  $f(X)$ to  $f[X(1-t)]$.

\begin{thm}\label{t J formula}
For any partition $\mu$ of $\ell$, the integral form Macdonald
polynomial $J_{\mu }(X;q,t)$ is given by
\begin{align}\label{et:J formula}
J_{\mu}(X;q,t) = t^{\mathsf{n}(\mu )} \mathbf{e}'_X\Big(z_1\cdots z_l
\prod_{\alpha_{ij} \in R_+} \!\! (1-z_i/z_j)\
\widehat{\HS}_{\mu}(\zz;q,t^{-1})\Big),
\end{align}
where
$\widehat{\HS}_\mu(\zz;q,t)$
is the
expression inside the $\sigmabold(\cdot)$ on the right side of
\eqref{et:mac formula a}.
Alternatively,
using informal
 notation similar to Remark~\ref{r:raising},
\begin{multline}\label{et:raising op}
J_{\mu}(X;q,t) = t^{\mathsf{n}(\mu )} \Bigg( \frac{\prod_{\alpha_{ij} \in
R_+}(1-\mathsf{R}_{ij}) \prod_{\alpha_{ij} \in \widehat{R}_{\mu}}
\big(1-q\spa t^{-1}\, \mathsf{R}_{ij} \big)} {\prod_{\alpha_{ij} \in
R_+} \big(1-q \, \mathsf{R}_{ij}\big)
\prod_{\alpha_{ij} \in R_\mu} \big(1-t^{-1} \, \mathsf{R}_{ij}\big)}\Bigg)  \\
 \times \prod_{\alpha_{ij} \in R_\mu \setminus \widehat{R}_\mu }
\big(1- q^{\arm(\bx{\mu}{i})+1}\spa t^{\leg(\bx{\mu}{i})}
\mathsf{R}_{ij} \big) \cdot e_{1^l}[X(1-t)],
\end{multline}
where $\mathsf{R}_{ij}$ acts on subscripts of $e_{\gamma}[X(1-t)]$ by
$\mathsf{R}_{ij} \gamma = \gamma + \epsilon_i - \epsilon_j$.
\end{thm}
\begin{proof}
By Theorem \ref{t mac formula} and \eqref{e:Weyl-symmetrizer 2},
\begin{align}
\Htild _{\mu }(X;q,t) = \omega \spa \mathbf{h}_X\big(z_1\cdots z_l
\prod_{\alpha_{ij} \in R_+} \!\! (1-z_i/z_j)\
\widehat{\HS}_{\mu}\big).
\end{align}
The result then follows from \eqref{e J in terms of Htild} and the
definition of $\mathbf{e}'_X$.
\end{proof}

We also record the consequences of
Corollary \ref{cor:mac formula b 2}
and
Corollary \ref{cor:m1 mac formula} for the integral form Macdonald
polynomials.
\begin{thm}
For any partition $\mu$ of $l$ and any rearrangement $\beta$ of
$\mu^*$,
\begin{align}
J_{\mu}(X;q,t) = t^{\mathsf{n}(\mu )} \mathbf{e}'_X\Big(z_1\cdots z_l
\prod_{\alpha_{ij} \in R_+} \!\! (1-z_i/z_j)\
\widehat{\HS}_{\beta ^{*}}(\zz;q,t^{-1})\Big),
\end{align}
where $\widehat{\HS}_{\beta ^{*}}(\zz;q,t)$ is the expression inside
$\sigmabold(\cdot)$ on the right side of \eqref{et:mac formula a 2}.
\end{thm}
\begin{thm}
Let $m \in \ZZ_{+}$, $\mu$ be a partition of $d$, and $l=dm$.  Then
for any rearrangement $\beta$ of $\mu^*$,
\begin{align}
J_{\mu}(X;q,t) = t^{m \spa \mathsf{n}(\mu)}
\mathbf{e}'_X\Big(\prod_{\alpha_{ij} \in R_+(\GL_l)} \!\! (1-z_i/z_j)\
\widehat{\HS}^{m,1}_{\beta ^{*}}(\zz;q,t^{-1})\Big),
\end{align}
where $\widehat{\HS}^{m,1}_{\beta ^{*}}(\zz;q,t)$ is the expression
inside $\sigmabold(\cdot)$ on the right side of \eqref{et:mac formula
mn a}, with $n$ set to 1.
\end{thm}


\begin{thebibliography}{10}

\bibitem{BeGaHaTe99}
F.~Bergeron, A.~M. Garsia, M.~Haiman, and G.~Tesler, \emph{Identities and
  positivity conjectures for some remarkable operators in the theory of
  symmetric functions}, Methods Appl. Anal. \textbf{6} (1999), no.~3, 363--420,
  Dedicated to Richard A.\ Askey on the occasion of his 65th birthday, Part
  III.

\bibitem{BeGaSeXi16}
Francois Bergeron, Adriano Garsia, Emily Sergel~Leven, and Guoce Xin,
  \emph{Compositional {$(km,kn)$}-shuffle conjectures}, Int. Math. Res. Not.
  IMRN (2016), no.~14, 4229--4270.

\bibitem{BHMPS-lw}
Jonah Blasiak, Mark Haiman, Jennifer Morse, Anna Pun, and George Seelinger,
  \emph{Dens, nests and the {L}oehr-{W}arrington conjecture}, 2021,
  \mbox{arXiv:2112.07070 [math.CO]}.

\bibitem{BHMPS-llt}
\bysame, \emph{{LLT} polynomials in the {S}chiffmann algebra}, 2021,
  \mbox{arXiv:2112.07063 [math.CO]}.

\bibitem{BHMPS-delta}
Jonah Blasiak, Mark Haiman, Jennifer Morse, Anna Pun, and George~H. Seelinger,
  \emph{A proof of the extended delta conjecture}, Forum Math. Pi \textbf{11}
  (2023), Paper No. e6, 28.

\bibitem{BHMPS-paths}
\bysame, \emph{A shuffle theorem for paths under any line}, Forum Math. Pi
  \textbf{11} (2023), Paper No. e5, 38.

\bibitem{BMPS2}
Jonah Blasiak, Jennifer Morse, Anna Pun, and Daniel Summers,
  \emph{{$k$}-{S}chur expansions of {C}atalan functions}, Adv. Math.
  \textbf{371} (2020), 107209, 39.

\bibitem{BroerNormality}
Bram Broer, \emph{Normality of some nilpotent varieties and cohomology of line
  bundles on the cotangent bundle of the flag variety}, Lie theory and
  geometry, Progr. Math., vol. 123, Birkh\"auser Boston, Boston, MA, 1994,
  pp.~1--19.

\bibitem{BurbSchi12}
Igor Burban and Olivier Schiffmann, \emph{On the {H}all algebra of an elliptic
  curve, {I}}, Duke Math. J. \textbf{161} (2012), no.~7, 1171--1231.

\bibitem{GarsHaim96}
A.~M. Garsia and M.~Haiman, \emph{Some natural bigraded {$S\sb n$}-modules and
  {$q,t$}-{K}ostka coefficients}, Electron. J. Combin. \textbf{3} (1996),
  no.~2, Research Paper 24, approx.\ 60 pp.\ (electronic), The Foata
  Festschrift.

\bibitem{HaHaLo05}
J.~Haglund, M.~Haiman, and N.~Loehr, \emph{A combinatorial formula for
  {M}acdonald polynomials}, J. Amer. Math. Soc. \textbf{18} (2005), no.~3,
  735--761 (electronic).

\bibitem{HaHaLo08}
\bysame, \emph{A combinatorial formula for nonsymmetric {M}acdonald
  polynomials}, Amer. J. Math. \textbf{130} (2008), no.~2, 359--383.

\bibitem{HaHaLoReUl05}
J.~Haglund, M.~Haiman, N.~Loehr, J.~B. Remmel, and A.~Ulyanov, \emph{A
  combinatorial formula for the character of the diagonal coinvariants}, Duke
  Math. J. \textbf{126} (2005), no.~2, 195--232.

\bibitem{Haiman01}
Mark Haiman, \emph{Hilbert schemes, polygraphs and the {M}acdonald positivity
  conjecture}, J. Amer. Math. Soc. \textbf{14} (2001), no.~4, 941--1006
  (electronic).

\bibitem{LaLeTh97}
Alain Lascoux, Bernard Leclerc, and Jean-Yves Thibon, \emph{Ribbon tableaux,
  {H}all-{L}ittlewood functions, quantum affine algebras, and unipotent
  varieties}, J. Math. Phys. \textbf{38} (1997), no.~2, 1041--1068.

\bibitem{LassalleSchlosser06}
Michel Lassalle and Michael Schlosser, \emph{Inversion of the {P}ieri formula
  for {M}acdonald polynomials}, Adv. Math. \textbf{202} (2006), no.~2,
  289--325.

\bibitem{Macdonald95}
I.~G. Macdonald, \emph{Symmetric functions and {H}all polynomials}, second ed.,
  The Clarendon Press, Oxford University Press, New York, 1995, With
  contributions by A.~Zelevinsky, Oxford Science Publications.

\bibitem{Mellit16}
Anton Mellit, \emph{Toric braids and $(m,n)$-parking functions}, 2016,
  \mbox{arXiv:1604.07456 [math.CO]}.

\bibitem{Milneclassicalpartitionfunctions}
S.~C. Milne, \emph{Classical partition functions and the {${\rm U}(n+1)$}
  {R}ogers-{S}elberg identity}, Discrete Math. \textbf{99} (1992), no.~1-3,
  199--246.

\bibitem{NoumShir12}
Masatoshi Noumi and Jun'ichi Shiraishi, \emph{A direct approach to the
  bispectral problem for the {R}uijsenaars-{M}acdonald q-difference operators},
  2012, \mbox{arXiv:1206.5364 [math.QA]}.

\bibitem{SchiVass13}
Olivier Schiffmann and Eric Vasserot, \emph{The elliptic {H}all algebra and the
  {$K$}-theory of the {H}ilbert scheme of {$\mathbb A^2$}}, Duke Math. J.
  \textbf{162} (2013), no.~2, 279--366.

\bibitem{ShimWeym00}
Mark Shimozono and Jerzy Weyman, \emph{Graded characters of modules supported
  in the closure of a nilpotent conjugacy class}, European J. Combin.
  \textbf{21} (2000), no.~2, 257--288.

\bibitem{Shiraishi05}
Jun'ichi Shiraishi, \emph{A conjecture about raising operators for {M}acdonald
  polynomials}, Lett. Math. Phys. \textbf{73} (2005), no.~1, 71--81.

\bibitem{Tamvakis11}
Harry Tamvakis, \emph{Giambelli, {P}ieri, and tableau formulas via raising
  operators}, J. Reine Angew. Math. \textbf{652} (2011), 207--244.

\bibitem{Thomas81}
Gl\^{a}nffrwd~P. Thomas, \emph{A note on {Y}oung's raising operator}, Canadian
  J. Math. \textbf{33} (1981), no.~1, 49--54.

\bibitem{Weyman89}
J.~Weyman, \emph{The equations of conjugacy classes of nilpotent matrices},
  Invent. Math. \textbf{98} (1989), no.~2, 229--245.

\end{thebibliography}

\providecommand{\bysame}{\leavevmode\hbox to3em{\hrulefill}\thinspace}

\end{document}